\documentclass[12pt]{amsart}
\usepackage{amssymb,euscript}
    \newtheorem{thm}{Theorem}                     [section]
    \newtheorem{prop}[thm]{Proposition}
    \newtheorem{lemma}[thm]{Lemma}
    \newtheorem{cor}[thm]{Corollary}
    \newtheorem{defn}[thm]{Definition}                 
    \newtheorem{rems}[thm]{Remark}                     

\newcommand{\ndef}{\newcommand*}
\def\rndef{\renewcommand}

\ndef{\myaddress}[1]{\begin{center} \it\tiny #1 \end{center}}

\ndef{\clA}{{\mathcal A}}
\ndef{\clB}{{\mathcal B}}
\ndef{\clE}{{\mathcal E}}
\ndef{\clH}{{\mathcal H}}
\ndef{\clK}{{\mathcal K}}
\ndef{\clL}{{\mathcal L}}
\ndef{\clN}{{\mathcal N}}
\ndef{\mbR}{{\mathbb R}}
\ndef{\mbC}{{\mathbb C}}
\ndef{\mbN}{{\mathbb N}}

\ndef{\eps}{\varepsilon}

\let\geq\geqslant
\let\leq\leqslant

\ndef{\lims}[1]{\lim\limits_{#1}}
\ndef{\sums}[1]{\sum\limits_{#1}}
\ndef{\ints}[1]{\int\limits_{#1}}
\ndef{\sups}[1]{\sup\limits_{#1}}
\ndef{\liminfty}[1]{\lims{#1\to\infty}}
\ndef{\suminf}[1]{\sums{#1=1}^\infty}

\ndef{\Tr}{\operatorname{Tr}}        

\ndef{\sflow}{\operatorname{sf}}     
\ndef{\sign}{\operatorname{sign}}    

\ndef{\rslv}[1]{R_z(#1)}      
\ndef{\HH}{H}                 
\ndef{\tHH}{\tilde \HH}       
\ndef{\VV}{V}                 
\ndef{\psif}[1]{#1^{[1]}} 
\ndef{\CPlus}[1]{W_{#1}(\mbR)}

\ndef{\hilb}{\clH}                     
\ndef{\hilbasargument}{(\hilb)} 
\ndef{\clBH}{\clB\hilbasargument}              

\ndef{\LpN}[1]{\clL^{#1}(\clN,\tau)}     
\ndef{\clAND}{(\clA,\clN,D)}             
\ndef{\clKN}{{\clK(\clN,\tau)}}          
\ndef{\oind}[2]{{\rm \tau\mbox{-}ind}_{#1\mbox{-}#2}}   
\ndef{\vNa}{von Neumann algebra}         
\ndef{\nsf}{faithful normal semifinite } 
\ndef{\taubrs}[1]{\tau\brackets{#1}}     
\ndef{\sqbrs}[1]{\left[#1\right]}        
\ndef{\IIinfty}{$\mathrm{II}_\infty$\ }
\ndef{\seque}[1]{\ensuremath{\{#1_j\}_{j=1}^\infty}}    
\ndef{\norm}[1]{\left\Vert#1\right\Vert}    
\ndef{\TrNorm}[1]{\norm{#1}_1}              
\ndef{\InftyNorm}[1]{\norm{#1}_\infty}      
\ndef{\abs}[1]{\left\lvert#1\right\rvert}   
\ndef{\set}[1]{\left\{#1\right\}}           
\ndef{\brackets}[1]{\left(#1\right)}        
\ndef{\brs}[1]{\brackets{#1}}               
\ndef{\scal}[2]{\left\la #1,#2\right\ra}               
\ndef{\csupp}{c}                           
\ndef{\mytimes}{\!\times\!}
\ndef{\then}{\Rightarrow}
\ndef{\la}{\langle}
\ndef{\ra}{\rangle}

\let\LatexCite=\cite  

\let\ifnumref\iftrue 

\ndef{\ifuncited}[4]{\expandafter\ifx\csname used#4\endcsname\relax}

\ndef{\ifcited}[4]{\expandafter\ifx\csname used#4\endcsname\relax\else}

  \ndef{\papertitle}[1]{ \emph{#1}, }
  \ndef{\paperauthor}[2]{#2}  
  \ndef{\pbbi}[9]{%
      \ifcited{#1}{#2}{#3}{#5}%
        \ifnumref%
          \bibitem{#5}\paperauthor{#1}{#6},\papertitle{#7}#8.%
        \else%
          \advance #9 by 1%
          \ifnum#9<1%
            \bibitem[#4]{#5}\paperauthor{#1}{#6}, \papertitle{#7}#8.%
          \else%
            \bibitem[#4$\!_{\the#9}\!$]{#5}\paperauthor{#1}{#6},\papertitle{#7}#8.%
          \fi%
        \fi%
      \fi%
  }
  \ndef{\mbbi}[8]{%
     \ifcited{#1}{#2}{#3}{#5}%
        \ifnumref%
          \bibitem{#5}\paperauthor{#1}{#6},\papertitle{#7}#8.%
        \else%
          \bibitem[#4]{#5}\paperauthor{#1}{#6},\papertitle{#7}#8.%
        \fi%
     \fi%
  }

\ndef{\AddCite}[1]{%
   \ifuncited{0}{0}{0}{#1}%
     \expandafter\gdef\csname used#1\endcsname {}%
   \fi%
}

\def\ProcessCite#1,{%
     \ifx\relax#1%
         \let\next=\relax%
     \else%
         \AddCite{#1}%
         \let\next=\ProcessCite%
     \fi%
     \next%
}

\ndef{\AddCites}[1]{\ProcessCite#1,\relax,}

\ndef{\CiteWithoutExtension}[1]{%
   \AddCites{#1}%
   \LatexCite{#1}%
}

\def\CiteWithExtension[#1]#2{%
   \AddCites{#2}%
   \LatexCite[#1]{#2}%
}

\ndef{\CleverCite}{%
    \ifx\NChar[ %
       \let\MyCite=\CiteWithExtension %
    \else %
       \let\MyCite=\CiteWithoutExtension %
    \fi %
    \MyCite%
}

\renewcommand{\cite}{\futurelet\NChar\CleverCite}

      \ndef{\volume}[1]{{\bf #1}}
      \ndef{\VolYearPP}[3]{\ifnum#2=0 (to appear)\else\volume{#1} (#2), #3\fi}
      \ndef{\VolNoYearPP}[4]{\ifnum#3=0 (to appear)\else\volume{#1} #2 (#3), #4\fi}
      \ndef{\libcode}[1]{}

\ndef{\jnActaMath}[3]{Acta Math. \VolYearPP{#1}{#2}{#3}}                       
\ndef{\jnAdvMath}[3]{Adv. in~Math. \VolYearPP{#1}{#2}{#3}}                     
\ndef{\jnAlgAnal}[3]{Algebra i~Analiz \VolYearPP{#1}{#2}{#3}}
\ndef{\jnAmerMathMonth}[3]{Amer. Math. Monthly \VolYearPP{#1}{#2}{#3}}         
\ndef{\jnAnnMath}[4]{Ann. of~Math. \VolNoYearPP{#1}{#2}{#3}{#4}}               
\ndef{\jnAnalMath}[3]{J. Anal. Math. \VolYearPP{#1}{#2}{#3}}                   
\ndef{\jnBullLondMathSoc}[3]{Bull. London Math. Soc. \VolYearPP{#1}{#2}{#3}}   
\ndef{\jnBullAMS}[3]{Bull. Amer. Math. Soc. \VolYearPP{#1}{#2}{#3}}   
\ndef{\jnCanMathBull}[3]{Canad. Math. Bull. \VolYearPP{#1}{#2}{#3}}            
\ndef{\jnCanMath}[4]{Canad. J.~Math. \VolNoYearPP{#1}{#2}{#3}{#4}}             
\ndef{\jnCommMathPhys}[3]{Comm. Math. Phys \VolYearPP{#1}{#2}{#3}}             
\ndef{\jnCommPDE}[3]{Comm. Partial Differential Equations \VolYearPP{#1}{#2}{#3}}             
\ndef{\jnComptRendue}[3]{C.\,R.~Acad. Sci. Paris S\'er. A-B \VolYearPP{#1}{#2}{#3}}      
\ndef{\jnDiffGeom}[3]{J.~Diff. Geom. \VolYearPP{#1}{#2}{#3}}                   
\ndef{\jnErgodicTheory}[3]{Ergodic Theory and Dynamical Systems \VolYearPP{#1}{#2}{#3}} 
\ndef{\jnFuncAnal}[3]{J.~Functional Analysis \VolYearPP{#1}{#2}{#3}}           
\ndef{\jnFunkAnalPril}[4]{Функциональный анализ и его приложения \VolNoYearPP{#1}{#2}{#3}{#4}}  
\ndef{\jnGAFA}[3]{GAFA \VolYearPP{#1}{#2}{#3}}                                 
\ndef{\jnIHES}[3]{IHES Publ. Math. (Paris) \VolYearPP{#1}{#2}{#3}}             
\ndef{\jnIEOT}[3]{Int. Equations and~Oper. Theory \VolYearPP{#1}{#2}{#3}} 
\ndef{\jnIsrMath}[3]{Israel J.~Math. \VolYearPP{#1}{#2}{#3}}                   
\ndef{\jnKTheory}[3]{K-Theory \VolYearPP{#1}{#2}{#3}}                          
\ndef{\jnLetMathPhys}[3]{Lett. Math. Phys. \VolYearPP{#1}{#2}{#3}}             
\ndef{\jnMathAnn}[3]{Math. Ann. \VolYearPP{#1}{#2}{#3}}                        
\ndef{\jnMathAnalAppl}[3]{J.~Math. Anal. and Appl. \VolYearPP{#1}{#2}{#3}}     
\ndef{\jnMathNachr}[3]{Math. Nachr. \VolYearPP{#1}{#2}{#3}}
\ndef{\jnMathPhys}[3]{J. Math. Phys. \VolYearPP{#1}{#2}{#3}}
\ndef{\jnOperTheory}[3]{J.~Operator Theory \VolYearPP{#1}{#2}{#3}}             
\ndef{\jnPacJMath}[3]{Pacific J.~Math. \VolYearPP{#1}{#2}{#3}}                  
\ndef{\jnPositivity}[3]{Positivity \VolYearPP{#1}{#2}{#3}}
\ndef{\jnProcAmerMS}[3]{Proc. Amer. Math. Soc. \VolYearPP{#1}{#2}{#3}}         
\ndef{\jnProcCambPhilSoc}[3]{Math. Proc. Camb. Phil. Soc. \VolYearPP{#1}{#2}{#3}}
\ndef{\jnReineAngew}[3]{J.~Reine Angew. Math. \VolYearPP{#1}{#2}{#3}}          
\ndef{\jnTokyoMath}[3]{Tokyo J.~Math. \VolYearPP{#1}{#2}{#3}}
\ndef{\jnTopology}[3]{Topology \VolYearPP{#1}{#2}{#3}}
\ndef{\jnTransAmerMathSoc}[3]{Trans. Amer. Math. Soc. \VolYearPP{#1}{#2}{#3}}
\ndef{\jnIzvANSSSR}[3]{Izv. Akad. Nauk SSSR, Ser. Mat. \VolYearPP{#1}{#2}{#3}}
\ndef{\jnIzvVyshUchZav}[3]{Izv. Vyssh. Uch. Zav., Mat. \VolYearPP{#1}{#2}{#3} (Russian)}
\ndef{\jnIzdatLenUniv}[2]{Izdat. Leningrad. Univ., Leningrad, (#1), #2 (Russian)}
\ndef{\jnFieldsInsComm}[3]{Fields Inst. Comm. \VolYearPP{#1}{#2}{#3}}
\ndef{\jnDoklANSSSR}[3]{Dokl. Akad. Nauk SSSR \VolYearPP{#1}{#2}{#3}}
\ndef{\jnMatZametki}[3]{Matem. zametki \VolYearPP{#1}{#2}{#3}}
\ndef{\jnRussMathSurvey}[3]{Russian Math. Surveys \VolYearPP{#1}{#2}{#3}}
\ndef{\jnSibMathJ}[3]{Sib. Math.~J. \VolYearPP{#1}{#2}{#3}}
\ndef{\jnSovMath}[3]{J.~Soviet math. \VolYearPP{#1}{#2}{#3}}
\ndef{\jnTransMoscMathSoc}[3]{Trans. Moscow Math. Soc. \VolYearPP{#1}{#2}{#3}}
\ndef{\jnUMN}[3]{UMN \VolYearPP{#1}{#2}{#3}}

\ndef{\bkTransMathMon}[2]{Trans. Math. Monographs, AMS, \volume{#1}, #2}

\ndef{\pbBirkhauser}[1]{Birkh\"auser, Boston, #1}
\ndef{\pbFactorial}[1]{Moscow, Factorial, #1}
\ndef{\pbGauthier}[1]{Gauthier-Villars, Paris, #1}
\ndef{\pbNauka}[1]{Moscow, Nauka, #1 (Russian)}
\ndef{\pbNaukaR}[1]{Москва, Наука, #1}
\ndef{\pbPrinceton}[1]{Princeton University Press, Princeton, New Jersey, #1}
\ndef{\pbPublPerish}[1]{Publish or Perish Inc., Berkeley, #1}
\ndef{\pbSpringer}[1]{Springer-Verlag, #1}

\ndef{\myauthor}[1]{\mbox{#1}}

\ndef{\Ahiezer}{\myauthor{N.\,I.\,Ahiezer }}
\ndef{\Arazy}{\myauthor{J.\,Arazy}}
\ndef{\Astashkin}{\myauthor{S.\,V.\,Astashkin}}
\ndef{\Atiyah}{\myauthor{M.\,Atiyah}}
\ndef{\Avron}{\myauthor{J.\,Avron}}
\ndef{\Azamov}{\myauthor{N.\,A.\,Azamov}}
\ndef{\Banach}{\myauthor{S.\,Banach}}
\ndef{\Benameur}{\myauthor{M-T.\,Benameur}}
\ndef{\Bennett}{\myauthor{C.\,Bennett}}
\ndef{\Berezin}{\myauthor{F.\,A.\,Berezin}}
\ndef{\Berline}{\myauthor{N.\,Berline}}
\ndef{\Birman}{\myauthor{M.\,Sh.\,Birman}}
\ndef{\Blackadar}{\myauthor{B.\,Blackadar}}
\ndef{\Bogolyubov}{\myauthor{N.\,N.\,Bogolyubov}}
\ndef{\Bonsall}{\myauthor{F.\,F.\,Bonsall}}
\ndef{\BoosBavnbek}{\myauthor{B.\,Boo$\beta$-Bavnbek}}
\ndef{\Bott}{\myauthor{R.\,Bott}}
\ndef{\Bratteli}{\myauthor{O.\,Bratteli}}
\ndef{\Bredon}{\myauthor{G.\,E.\,Bredon}}
\ndef{\Breuer}{\myauthor{M.\,Breuer}}
\ndef{\Brown}{\myauthor{L.\,G.\,Brown}}
\ndef{\Bruneau}{\myauthor{V.\,Bruneau}}
\ndef{\Buslaev}{\myauthor{V.\,S.\,Buslaev}}
\ndef{\Carey}{\myauthor{A.\,L.\,Carey}}
\ndef{\CareyRW}{\myauthor{R.\,W.\,Carey}} 
\ndef{\Cartan}{\myauthor{H.\,Cartan}}
\ndef{\Chilin}{\myauthor{V.\,I.\,Chilin}}
\ndef{\Coburn}{\myauthor{L.\,A.\,Coburn}}
\ndef{\Connes}{\myauthor{A.\,Connes}}
\ndef{\Cornfeld}{\myauthor{I.\,P.\,Cornfeld}}
\ndef{\Daletskii}{\myauthor{Yu.\,L.\,Daletski\u\i}}   
\ndef{\Dixmier}{\myauthor{J.\,Dixmier}}
\ndef{\DoddsPG}{\myauthor{P.\,G.\,Dodds}}
\ndef{\DoddsTK}{\myauthor{T.\,K.\,Dodds}}
\ndef{\Douglas}{\myauthor{R.\,G.\,Douglas}}
\ndef{\Dubrovin}{\myauthor{B.\,A.\,Dubrovin}}
\ndef{\Dugundji}{\myauthor{J.\,Dugundji}}
\ndef{\Duncan}{\myauthor{J.\,Duncan}}
\ndef{\Dunford}{\myauthor{N.\,Dunford}}
\ndef{\Dykema}{\myauthor{K.\,J.\,Dykema}}
\ndef{\Edwards}{\myauthor{R.\,E.\,Edwards}}
\ndef{\Eilenberg}{\myauthor{S.\,Eilenberg}}
\ndef{\Fack}{\myauthor{T.\,Fack}} 
\ndef{\Faddeev}{\myauthor{L.\,D.\,Faddeev}}
\ndef{\Farber}{\myauthor{M.\,Farber}}
\ndef{\Farforovskaya}{\myauthor{Yu.\,B.\,Farforovskaya}}
\ndef{\Federer}{\myauthor{H.\,Federer}}
\ndef{\Fedosov}{\myauthor{B.\,V.\,Fedosov}}
\ndef{\Figiel}{\myauthor{T.\,Figiel}} 
\ndef{\Figueroa}{\myauthor{H.\,Figueroa}}
\ndef{\Fillmore}{\myauthor{P.\,A.\,Fillmore}}
\ndef{\Fomenko}{\myauthor{A.\,T.\,Fomenko}} 
\ndef{\Fomin}{\myauthor{S.\,V.\,Fomin}}
\ndef{\Frohlich}{\myauthor{J.\,Fr\"ohlich}}
\ndef{\Fuglede}{\myauthor{B.\,Fuglede}}
\ndef{\Furutani}{\myauthor{K.\,Furutani}}
\ndef{\Gelfand}{\myauthor{I.\,M.\,Gelfand}}
\ndef{\Gesztesy}{\myauthor{F.\,Gesztesy}}     
\ndef{\Getzler}{\myauthor{E.\,Getzler}} 
\ndef{\Gilkey}{\myauthor{P.\,B.\,Gilkey}}
\ndef{\Gitler}{\myauthor{S.\,Gitler}}
\ndef{\Glazman}{\myauthor{I.\,M.\,Glazman}}
\ndef{\Glimm}{\myauthor{J.\,Glimm}}
\ndef{\Gohberg}{\myauthor{I.\,C.\,Gohberg}}
\ndef{\Golze}{\myauthor{F.\,Golze}}
\ndef{\GraciaBondia}{\myauthor{J.\,M.\,Gracia-Bond\'{i}a}}
\ndef{\Greenleaf}{\myauthor{F.\,P.\,Greenleaf}}
\ndef{\Gromov}{\myauthor{M.\,Gromov}}
\ndef{\Gunning}{\myauthor{R.\,C.\,Gunning}}
\ndef{\Haagerup}{\myauthor{U.\,Haagerup}}
\ndef{\Haag}{\myauthor{R.\,Haag}}
\ndef{\Halmos}{\myauthor{Halmos}}
\ndef{\Hardy}{\myauthor{G.\,H.\,Hardy}}
\ndef{\Higson}{\myauthor{N.\,Higson}}  
\ndef{\Hoermander}{\myauthor{L.\,Hoermander}} 
\ndef{\Hoffman}{\myauthor{K.\,Hoffman}} 
\ndef{\Ito}{\myauthor{K.\,Ito}}
\ndef{\Jaffe}{\myauthor{A.\,Jaffe}}
\ndef{\James}{\myauthor{I.\,M.\,James}}
\ndef{\Javrjan}{\myauthor{V.\,A.\,Javrjan}}
\ndef{\Kadison}{\myauthor{R.\,V.\,Kadison}}
\ndef{\Kalton}{\myauthor{N.\,J.\,Kalton}} 
\ndef{\Kato}{\myauthor{T.\,Kato}} 
\ndef{\Kobayashi}{\myauthor{S.\,Kobayashi}}
\ndef{\Koplienko}{\myauthor{L.\,S.\,Koplienko}}
\ndef{\Korotyaev}{\myauthor{E.\,Korotyaev}}
\ndef{\Kosaki}{\myauthor{H.\,Kosaki}}
\ndef{\KreinMG}{\myauthor{M.\,G.\,Kre\u\i n}}
\ndef{\KreinSG}{\myauthor{S.\,G.\,Kre\u\i n}}
\ndef{\Leichtnam}{\myauthor{E.\,Leichtnam}}
\ndef{\Lesch}{\myauthor{M.\,Lesch}}
\ndef{\Lesniewski}{\myauthor{A.\,Lesniewski}}
\ndef{\Levitan}{\myauthor{B.\,M.\,Levitan}}
\ndef{\Lidskii}{\myauthor{V.\,B.\,Lidskii}}
\ndef{\Lifshits}{\myauthor{I.\,M.\,Lifshits}}
\ndef{\Lindenstrauss}{\myauthor{J.\,Lindenstrauss}}
\ndef{\Loday}{\myauthor{J.-L.\,Loday}}
\ndef{\Lord}{\myauthor{S.\,Lord}}      
\ndef{\Lorentz}{\myauthor{G.\,Lorentz}}
\ndef{\Magnus}{\myauthor{W.\,Magnus}}
\ndef{\Makarov}{\myauthor{K.\,A.\,Makarov}}
\ndef{\Mathai}{\myauthor{V.\,Mathai}}         
\ndef{\McKean}{\myauthor{H.\,P.\,McKean}}
\ndef{\Mishchenko}{\myauthor{A.\,S.\,Mishchenko}}
\ndef{\Moore}{\myauthor{C.\,C.\,Moore}}
\ndef{\Moscovici}{\myauthor{H.\,Moscovici}}  
\ndef{\Motovilov}{\myauthor{A.\,K.\,Motovilov}}
\ndef{\Moyer}{\myauthor{R.\,D.\,Moyer}}
\ndef{\Naboko}{\myauthor{S.\,N.\,Naboko}}
\ndef{\Narasimhan}{\myauthor{R.\,Narasimhan}}
\ndef{\Nomizu}{\myauthor{K.\,Nomizu}}
\ndef{\Novikov}{\myauthor{S.\,P.\,Novikov}}
\ndef{\Osterwalder}{\myauthor{K.\,Osterwalder}}
\ndef{\Patodi}{\myauthor{V.\,Patodi}}
\ndef{\Pagter}{\myauthor{B.\,de~Pagter}}  
\ndef{\Pavlov}{\myauthor{B.\,S.\,Pavlov}}
\ndef{\Pedersen}{\myauthor{G.\,K.\,Pedersen}}
\ndef{\Peller}{\myauthor{V.\,V.\,Peller}}
\ndef{\Perera}{\myauthor{V.\,S.\,Perera}}
\ndef{\Petunin}{\myauthor{Ju.\,I.\,Petunin}}
\ndef{\Phillips}{\myauthor{J.\,Phillips}}  
\ndef{\Piazza}{\myauthor{P.\,Piazza}}   
\ndef{\Pincus}{\myauthor{J.\,D.\,Pincus}}   
\ndef{\Poincare}{Poincar\'e}
\ndef{\Postnikov}{\myauthor{M.\,M.\,Postnikov}} 
\ndef{\Prinzis}{\myauthor{R.\,Prinzis}}
\ndef{\Privalov}{\myauthor{I.\,I.\,Privalov}}
\ndef{\Pushnitski}{\myauthor{A.\,B.\,Pushnitski}} 
\ndef{\Raeburn}{\myauthor{I.\,Raeburn}}
\ndef{\Raikov}{\myauthor{G.\,Raikov}}
\ndef{\Reed}{\myauthor{M.\,Reed}}
\ndef{\Rennie}{\myauthor{A.\,Rennie}}
\ndef{\Rickart}{\myauthor{C.\,E.\,Rickart}}
\ndef{\Riesz}{\myauthor{F.\,Riesz}}
\ndef{\Ringrose}{\myauthor{J.\,Ringrose}}
\ndef{\Robinson}{\myauthor{D.\,Robinson}}
\ndef{\Rossi}{\myauthor{H.\,Rossi}}
\ndef{\Rudin}{\myauthor{W.\,Rudin}}
\ndef{\Ruelle}{\myauthor{D.\,Ruelle}}
\ndef{\Ruzhansky}{\myauthor{M.\,Ruzhansky}}
\ndef{\Sakai}{\myauthor{Sh.\,Sakai}}
\ndef{\Sargsjan}{\myauthor{I.\,S.\,Sargsjan}}
\ndef{\Sato}{\myauthor{H.\,Sato}}
\ndef{\Schaeffer}{\myauthor{D.\,G.\,Schaeffer}}
\ndef{\Schluchtermann}{\myauthor{G.\,Schluchtermann}}
\ndef{\Schochet}{\myauthor{C.\,Schochet}}
\ndef{\Schrodinger}{\myauthor{E.\,Schr\"odinger}}
\ndef{\Schrohe}{\myauthor{E.\,Schrohe}}
\ndef{\Schwartz}{\myauthor{J.\,T.\,Schwartz}}
\ndef{\Sedaev}{\myauthor{A.\,A.\,Sedaev}}
\ndef{\Seiler}{\myauthor{R.\,Seiler}}
\ndef{\Semenov}{\myauthor{E.\,M.\,Semenov}}
\ndef{\Shabat}{\myauthor{B.\,V.\,Shabat}}
\ndef{\Shafarevich}{\myauthor{I.\,R.\,Shafarevich}}
\ndef{\Sharpley}{\myauthor{R.\,Sharpley}}
\ndef{\Shilov}{\myauthor{G.\,E.\,Shilov}}
\ndef{\Shirkov}{\myauthor{D.\,V.\,Shirkov}}
\ndef{\Shubin}{\myauthor{M.\,A.\,Shubin}}
\ndef{\Silverman}{\myauthor{H.\,Silverman}}
\ndef{\Simon}{\myauthor{B.\,Simon}}
\ndef{\Sinai}{\myauthor{Ya.\,G.\,Sinai}}
\ndef{\Singer}{\myauthor{I.\,M.\,Singer}}
\ndef{\Solomyak}{\myauthor{M.\,Z.\,Solomyak}}
\ndef{\Soloviev}{\myauthor{Yu.\,P.\,Soloviev}}
\ndef{\Spivak}{\myauthor{M.\,Spivak}}
\ndef{\Stenkin}{\myauthor{V.\,V.\,Sten'kin}}
\ndef{\Stratila}{\myauthor{S.\,Stratila}}
\ndef{\Sucheston}{\myauthor{L.\,Sucheston}}
\ndef{\Sukochev}{\myauthor{F.\,A.\,Sukochev}}
\ndef{\Switzer}{\myauthor{R.\,M.\,Switzer}}
\ndef{\SzNagy}{\myauthor{B.\,Sz.-Nagy}}
\ndef{\Takesaki}{\myauthor{M.\,Takesaki}}
\ndef{\Taylor}{\myauthor{M.\,E.\,Taylor}}
\ndef{\Treves}{\myauthor{F.\,Treves}}
\ndef{\Troitsky}{\myauthor{E.\,V.\,Troitsky}}
\ndef{\Tzafriri}{\myauthor{L.\,Tzafriri}}
\ndef{\Varilly}{\myauthor{J.\,C.\,V\'{a}rilly}}
\ndef{\Vergne}{\myauthor{M.\,Vergne}}
\ndef{\Vladimirov}{\myauthor{V.\,S.\,Vladimirov}}
\ndef{\Voiculescu}{\myauthor{D.\,Voiculescu}}
\ndef{\Weiss}{\myauthor{G.\,Weiss}}
\ndef{\Wells}{\myauthor{R.\,O.\,Wells}}
\ndef{\Williams}{\myauthor{J.\,P.\,Williams}}
\ndef{\Winkler}{\myauthor{S.\,Winkler}}
\ndef{\Witten}{\myauthor{E.\,Witten}}
\ndef{\Wodzicki}{\myauthor{M.\,Wodzicki}}
\ndef{\Wojciechowski}{\myauthor{K.\,P.\,Wojciechowski}}
\ndef{\Yafaev}{\myauthor{D.\,R.\,Yafaev}}
\ndef{\Yosida}{\myauthor{K.\,Yosida}}
\ndef{\Zsido}{\myauthor{L.\,Zsido}}
\begin{document}
\begin{center} OPERATOR INTEGRALS, SPECTRAL SHIFT
   \\  AND~SPECTRAL~FLOW
\end{center}
\vskip 1 cm

\begin{center}
   \Azamov$^a$, \ \Carey$^{b,}$\footnote{Research partially supported by the Australian
Research Council \hfil\break \indent 2000
{\it Mathematics Subject Classification}: Primary 47A56, 47B49; Secondary 47A55, 46L51}, \ \DoddsPG$^{c,1}$, \ \Sukochev$^{d,1}$
\end{center}
\vskip 1 cm

\myaddress{$^{a,c,d}$ School of Informatics and Engineering,
   \\ Flinders University of South Australia,
   \\ Bedford Park, 5042, SA, Australia.
   \\ Email: azam0001,peter,sukochev@infoeng.flinders.edu.au}
\myaddress{$^{b}$ Mathematical Sciences Institute,
   \\ Australian National University,
   \\ Canberra, ACT 0200, Australia,
   \\ Email: acarey@maths.anu.edu.au}

\ndef{\intint}{\int\!\!\!\!\int}   
\ndef{\Pin}[1]{\Pi^{(#1)}}         
\ndef{\nunf}[1]{\nu^{(#1)}_f}      
\ndef{\ww}{\mbox{$*$-}~}           
\ndef{\LsInf}{\clL^{so^*}_\infty}
\rndef{\HH}{B}                     
\rndef{\VV}{V}
\rndef{\tHH}{{\widetilde\HH}}                
\ndef{\rr}{r}                      
\ndef{\epspi}{\sqrt{\frac{\eps}{\pi}}}
\rndef{\clAND}{(\clA,\clN,\HH)}

\vskip 1 cm
{\small
Abstract.

We present a new and simple approach to the theory of multiple
operator integrals that applies to unbounded operators affiliated with general \vNa s.
For semifinite \vNa s we give applications
to the Fr\'echet differentiation of operator functions that sharpen existing results,
and establish the Birman-Solomyak representation of the spectral shift function of M.\,G.\,Krein
in terms of an average of spectral measures in the type II setting.
We also exhibit a surprising connection between the spectral shift
function and spectral flow.
}
\vskip 1 cm

\section{Introduction}

In the seminal paper of Yu.\,L.\,Daletskii  and S.\,G.\,Krein \cite{DalKr},
the theory of multiple operator integrals emerged as an important tool in the
differentiation theory of operator functions and in perturbation theory.
On the other hand, an important concept in the theory of perturbations is
the spectral shift function which first arose in the work of \Lifshits\ \cite{Li52UMN}
in solid state theory and put on a firm mathematical basis by M.\,G.\,Krein \cite{Kr53MS}.

An important connection between these two theories was made by
\Birman\ and \Solomyak\ \cite{BS72SM} who showed that the theory of double operator integrals
led naturally to a new representation for the spectral shift function as an average of spectral measures.

A principal aim of this paper is to present a new approach to the theory of multiple operator integrals,
which provides a coherent path to the theory of differentiation of operator functions, the spectral shift function
and the theory of spectral flow in the setting of type II \vNa s.
Our approach is conceptually simpler
than those of \cite{dPS04FA,dPSW02FA}, although it applies to a somewhat
narrower class of functions.
On the other hand, since our approach does not depend on the vector-valued
integration theory against a finitely additive measure as in \cite{BS73DOI,dPS04FA,dPSW02FA},
it is also suitable for general (non-semifinite) \vNa s (see Section~\ref{S: MOI section}).
Our approach to the theory of multiple operator integrals is quite different from earlier
approaches to be found in \cite{DalKr,Pa,SolSt,St}.

During the final stages in the preparation of this paper, the authors became aware
of a preprint of \Peller\ \cite{Pel2}, where a similar approach to the theory of
multiple operator integrals is presented in the setting of type I \vNa s. While
Peller's approach in the type I case applies to the class of integral projective
tensor products, the present paper in the more general type II setting restricts attention
to the natural Wiener classes and thus permits us to show Fr\'echet differentiability
rather than Gateaux differentiability.
In particular, we strengthen the differentiation results of \cite{dPS04FA} by showing that
Gateaux differentiability can be replaced by Fr\'echet differentiability and we
show the existence of higher order Fr\'echet derivatives
of operator-valued functions.

Our present approach
also allows us to consider perturbations of self-adjoint operators which are affiliated with semifinite
\vNa s. This is a substantial difference with
\cite{dPS04FA,dPSW02FA}, which treated the more special case of (so-called) $\tau$-measurable
operators \cite{FK86PJM}. The necessity of avoiding
the latter restriction is especially clear in the applications.
An important ingredient is the recent extension to the type II setting of the Krein
spectral shift function \cite{ADS}. When combined with our development of multiple
operator integration together with the ideas of \cite{BS72SM}, we establish
a (type II) extension of the important Birman-Solomyak formula,
concerning spectral averaging (see Section~\ref{S: BirSol in vNa}).

Perhaps the most surprising connection of our study is
with the theory of spectral flow, which is
presented in Section~\ref{S: SF and SSF}. While the theory of spectral shift function (see the lectures
\cite{Kr63} and the survey \cite{BP98IEOT}) is a part of operator theory, the theory of spectral flow,
which originated in the work of \Atiyah, \Patodi\ and \Singer\ \cite{APS76} on a
generalization of the Atiyah-Singer index theorem, finds its proper
analytic setting in the framework of non-commutative geometry
created by \Connes\ \cite{CoNG}. One of the main results of the latter theory is
the odd local index theorem of \Connes\ and \Moscovici\ \cite{CM95GAFA} which
 has recently been developed in the type II setting \cite{CPRS2,CPRS3}.
This latter work inspired our
result (in Section~\ref{S: SF and SSF}) that the theory of the spectral shift
function and that of spectral flow coincide in the case of trace class perturbations.

\section{Notations and Preliminaries}
We denote by $\hilb$ a separable complex Hilbert space,
by $\clN$ a von Neumann algebra acting on $\hilb,$
by $\clBH$ the algebra of all bounded linear operators acting on $\hilb$
and by $\Tr$ the standard trace on $\hilb.$
    In case when $\clN$ is semifinite, we denote
by $\tau$ a \nsf trace on $\clN,$
by $L^1(\clN,\tau)$ the set of $\tau$-trace class operators affiliated with $\clN,$
by $\LpN{1}=L^1(\clN,\tau)\cap \clN$ the set of all bounded $\tau$-trace class operators
and by $\clKN$ the set of all $\tau$-compact operators (see \cite{FK86PJM}) from $\clN.$
    If $S$ is a measure space we denote
by $B(S)$ the set of all bounded measurable complex-valued functions on $S.$
The $so$-topology and $so^*$-topology denote, respectively, the strong operator topology and the strong$^*$ operator topology.
We denote the uniform norm on $\clBH$ by $\norm{\cdot}.$

If $\clE$ is a $*$-ideal in a \vNa\ $\clN$ which is complete in some norm $\norm{\cdot}_\clE,$ then we will call $\clE$
an invariant operator ideal (see \cite[Definition 1.8]{CPS03AdvM}) if
\\ \quad (1) $\norm{S}_\clE \geq \norm{S}$ for all $S \in \clE$,
\\ \quad (2) $\norm{S^*}_\clE = \norm{S}_\clE$ for all $S \in \clE$,
\\ \quad (3) $\norm{ASB}_\clE \leq \norm{A}\norm{S}_\clE\norm{B}$ for all $S \in \clE$ and $A,B \in \clN$ .

We say that an operator ideal $\clE$ has property (F) if, for all nets $\set{A_\alpha} \subset \clE$
such that there exist $A \in \clN$ for which $A_\alpha \to A$ in the $so^*$-topology
and $\norm{A_\alpha}_\clE \leq 1$ for all $\alpha,$ it follows that $A \in \clE$ and $\norm{A}_\clE \leq 1.$

If $\clE = \clN \cap E(\clN,\tau)$ for some rearrangement invariant Banach function space $E$ (see \cite{DDPS}) with
the Fatou property (that is, if \, $ 0 \leq f_\alpha \uparrow$ \, is an increasing net in $E$, $\sup\limits_\alpha\norm{f_\alpha}_E < \infty$
then $\sup\limits_\alpha f_\alpha$ exists in $E$ and $\norm{f_\alpha}_E \uparrow \norm{f}_E$), then
\cite[Proposition 1.6]{DDPS} together with Lemma~\ref{L: net VA(alpha) converges to VA} below shows that $\clE$ has the property (F).

Every \vNa\ with the uniform norm is an invariant operator ideal with property (F). If $\clN$ is a semifinite \vNa\ with a \nsf trace $\tau$
then the spaces $\clL^p(\clN,\tau), \ \clL^{p,+\infty}(\clN,\tau)$ (see \cite{KPS,DDP}) are invariant operator ideals with the property (F).

For any $C^1$-function $f\colon \mbR \to \mbC,$
we denote by $f^{[1]}$ the continuous function
$$f^{[1]}(\lambda_0,\lambda_1) = \frac {f(\lambda_1) - f(\lambda_0)}{\lambda_1-\lambda_0},$$ 
and for any $C^{n+1}$-function $f\colon \mbR \to \mbC$
$$
  f^{[n+1]}(\lambda_0,\ldots,\lambda_{n+1})
    = \frac {f^{[n]}(\lambda_0,\ldots,\lambda_{n-1},\lambda_{n+1})-f^{[n]}(\lambda_0,\ldots,\lambda_{n-1},\lambda_{n})}{\lambda_{n+1}-\lambda_{n}}.
$$
It is well known that $f^{[n]}$ is a symmetric function.

We denote by $\CPlus{n}$ the set of functions $f\in C^n(\mbR),$  such that the $j$-th derivative $f^{(j)}, \ j=0,\ldots,n,$
is the Fourier transform of a finite measure $m_f$ on $\mbR.$

The next lemma introduces a finite measure space which will be used in our
definition of multiple operator integrals in Section~\ref{S: MOI section} below.
\begin{lemma}\label{L: (Pi,nu) is a finite measure space}
If
\begin{gather*}
  \Pin{n} = \{(s_0,s_1,\ldots,s_n) \in \mbR^{n+1}\colon \abs{s_n}\leq\ldots\leq\abs{s_1}\leq\abs{s_0},
   \\ \sign(s_0) = \ldots = \sign(s_n)\},
\end{gather*}
and if  $f\in \CPlus{n},$ $\nunf{n}(s_0,\ldots,s_n) = \frac {i^n}{\sqrt{2\pi}} \,m_f(ds_0)\,\ldots\,ds_n,$
then $\brs{\Pin{n}, \nunf{n}}$ is a finite measure space.
\end{lemma}
\begin{proof} The total variation of the measure $\nunf{n}$ on the set $\Pin{n}$ (up to a constant) is equal to
\begin{align*}
    \int_{\Pin{n}} \abs{m_f(ds_0)}\,ds_1\, & \ldots\,ds_n
     = \int_{\mbR} \Delta_{s_0} \abs{m_f(ds_0)}
    \\ & = \frac 1{n!} \int_{\mbR} s_0^n \abs{m_f(ds_0)}
    \\ & = \frac 1{n!} \int_{\mbR} \abs{m_{f^{(n)}}(ds_0)} = \frac 1{n!} \norm{m_{f^{(n)}}},
\end{align*}
where $\Delta_{s_0}$ is the volume of the $n$-dimensional simplex of size $s_0.$
\end{proof}

We write for simplicity $\Pi = \Pin{1}$ and $\nu_f = \nunf{1}.$

The next two lemmas provide concrete representations (see section~\ref{S: MOI section})
for divided differences $f^{[n]}$ of functions belonging to the class $\CPlus{n}.$

\begin{lemma}\label{L: BS representation for f(x) - f(y)} If $f \in \CPlus 1,$ then
  \begin{gather}\label{F: psif = double}
    \psif{f}(\lambda_0,\lambda_1)
      = \intint_\Pi \alpha_0(\lambda_0,\sigma)\alpha_1(\lambda_1,\sigma)\,d\nu_f(\sigma),
  \end{gather}
   where $\sigma = (s_0,s_1),$ $\alpha_0(\lambda_0,\sigma) = e^{i(s_0-s_1)\lambda_0}, \quad \alpha_1(\lambda_1,\sigma) = e^{is_1\lambda_1},$
   \ $s_0, s_1 \in \mbR.$
\end{lemma}
\begin{proof}
  We have
  \begin{align*}
    \intint_\Pi \alpha_0(\lambda_0,\sigma) & \alpha_1(\lambda_1,\sigma)\,d\nu_f(\sigma)
    \\ & = \frac {i}{\sqrt{2\pi}} \int_\mbR m_f(ds_0) \int_0^{s_0} e^{is_0\lambda_0-is_1\lambda_0+is_1\lambda_1}\,ds_1
    \\ & = \frac {1}{(\lambda_0-\lambda_1)\sqrt{2\pi}} \int_\mbR m_f(ds_0)(e^{is_0\lambda_0} - e^{is_0\lambda_1})
    \\ & = \frac {1}{\lambda_0-\lambda_1} (f(\lambda_0) - f(\lambda_1)) = \psif{f}(\lambda_0,\lambda_1),
  \end{align*}
  where the repeated integral can be replaced
    by the double integral by Fubini's theorem and Lemma~\ref{L: (Pi,nu) is a finite measure space}.
\end{proof}

\begin{lemma}\label{L: integral formula for n-th divided difference} If $f\in \CPlus{n},$ then,
 for all $\lambda_0,\ldots,\lambda_n \in \mbR,$
  \begin{multline*}
    f^{[n]}(\lambda_0,\ldots,\lambda_n)
       \\ = \int_{\Pin{n}} e^{i((s_0-s_1)\lambda_0 + \ldots+(s_{n-1}-s_n)\lambda_{n-1}+s_n\lambda_n)}\,d\nunf{n}(s_0,\ldots,s_n).
  \end{multline*}
\end{lemma}
\begin{proof} By Lemma~\ref{L: BS representation for f(x) - f(y)} and induction, we have
  \begin{multline*}
  \int_{\Pin{n+1}} e^{i((s_0-s_1)\lambda_0+\ldots+(s_{n}-s_{n+1})\lambda_{n}+s_{n+1}\lambda_{n+1})} \,d\nunf{n+1}(s_0,\ldots,s_{n+1})
    \\  = \int_{\Pin{n}} e^{i((s_0-s_1)\lambda_0+\ldots+s_{n}\lambda_{n})}
    \brs{\int_0^{s_n} ie^{is_{n+1}(\lambda_{n+1} - \lambda_{n})}\,ds_{n+1}}\,d\nunf{n}(s_0,\ldots,s_{n})
    \\  = \frac 1{\lambda_{n+1}-\lambda_n} \int_{\Pin{n}} e^{i((s_0-s_1)\lambda_0+\ldots+s_{n}\lambda_{n})}
    \brs{e^{is_{n}(\lambda_{n+1} - \lambda_{n})}-1}\,d\nunf{n}(s_0,\ldots,s_{n})
    \\  = \frac 1{\lambda_{n+1}-\lambda_n} \brs{f^{[n]}(\lambda_0,\ldots,\lambda_{n-1},\lambda_{n+1}) - f^{[n]}(\lambda_0,\ldots,\lambda_{n-1},\lambda_{n})}
    \\  = f^{[n+1]}(\lambda_0,\ldots,\lambda_{n+1}).
  \end{multline*}
\end{proof}
\begin{lemma}\label{L: calculus exercise} If $f \in \CPlus {n+1},$ then,
 for all $\lambda_0,\ldots,\lambda_{n+1} \in \mbR,$
  \begin{multline*}
    f^{[n+1]}(\lambda_0,\ldots,\lambda_{n+1})
      \\ = i\int_{\Pin{n}}\int_0^{s_{j}-s_{j+1}}  e^{i((s_0-s_1)\lambda_0 + \ldots
           + u\lambda_{n+1} + (s_{j}-s_{j+1}-u)\lambda_{j}+ \ldots + s_n\lambda_n}
       \\  \times  \,du\,d\nunf{n}(s_0,\ldots,s_n).
  \end{multline*}
\end{lemma}
\begin{proof}
  The right hand side is equal to
     \begin{gather*}
       i\int_{\Pin{n}}
             e^{i((s_0-s_1)\lambda_0 + \ldots + (s_{j}-s_{j+1})\lambda_{j}+ \ldots + s_n\lambda_n)}
            \int_0^{s_{j}-s_{j+1}} e^{iu(\lambda_{n+1}-\lambda_{j})}
            \\ \times \,du\,d\nunf{n}(s_0,\ldots,s_n)
           \\ = \frac{1}{\lambda_{n+1}-\lambda_{j}} \int_{\Pin{n}}
           e^{i((s_0-s_1)\lambda_0 + \ldots + (s_{j}-s_{j+1})\lambda_{j}+ \ldots + s_n\lambda_n)}
           \\ \brs{e^{(s_{j}-s_{j+1})(\lambda_{n+1}-\lambda_{j})}-1} \,d\nunf{n}(s_0,\ldots,s_n)
           \\ = \frac{1}{\lambda_{n+1}-\lambda_{j}}
                \int_{\Pin{n}} \Big( e^{i((s_0-s_1)\lambda_0 + \ldots + (s_{j}-s_{j+1})\lambda_{n+1}+ \ldots + s_n\lambda_n)}
            \\    -  e^{i((s_0-s_1)\lambda_0 + \ldots + (s_{j}-s_{j+1})\lambda_{j}+ \ldots + s_n\lambda_n)} \Big)
               \,d\nunf{n}(s_0,\ldots,s_n)
           \\ = \frac{1}{\lambda_{n+1}-\lambda_{j}} \big( f^{[n]}(\lambda_0,\ldots,\lambda_{j-1},\lambda_{n+1},\lambda_{j+1},\ldots,\lambda_n)
               \\ - f^{[n]}(\lambda_0,\ldots,\lambda_{j-1},\lambda_{j},\lambda_{j+1},\ldots,\lambda_n) \big)
              = f^{[n+1]}(\lambda_0,\ldots,\lambda_{n+1}).
     \end{gather*}
\end{proof}
\begin{lemma}\label{L: net VA(alpha) converges to VA} Let $(\clN,\tau)$ be a semifinite \vNa.
If $A_\alpha \in \clN, \ \alpha \in I,$ is a uniformly bounded net
converging in the $so$-topology to an
operator $A \in \clN$ and if $\VV \in \LpN{1},$ then
the net $\set{A_\alpha\VV}_{\alpha\in I}$ 
converges to $A\VV$ 
in $L^{1}(\clN,\tau).$
\end{lemma}
\begin{proof} 
Without loss of generality, we can assume that $A=0.$
Since the net $\set{A_\alpha}_{\alpha \in I}$ is uniformly bounded, we have $A_\alpha \to 0$ in the
$\sigma$-strong operator topology (see e.g. \cite[Proposition 2.4.1]{BR1}). Since the $\sigma$-strong topology does not
depend on representation \cite[Theorem 2.4.23]{BR1}, it can be assumed that $\clN$ acts
on $L^{2}(\clN,\tau)$ in the left regular representation,
in particular $\norm{A_\alpha y}_2 \to 0$ for every $y \in \LpN{2}.$
Without loss of generality, we may suppose that $\VV \geq 0.$
Let $y = \VV^{1/2} \in \LpN{2}.$ Then
$$
  \taubrs{\abs{A_\alpha \VV}} = \taubrs{u_\alpha A_\alpha y^2}
  = \taubrs{A_\alpha y (u_\alpha^* y)^* } \leq \norm{A_\alpha y}_2 \cdot \norm{u_\alpha^* y}_2 \to 0,
$$
where $u^*_\alpha$ is the partial isometry from the polar decomposition of $A_\alpha \VV.$
\end{proof}
\begin{lemma}\label{L: the measure tau(AE(a,b)B) is count add} Let $A,B \in \clN$ and suppose
that one of these operators is $\tau$-trace-class. If $T=T^*$ is affiliated with $\clN$
and if $T = \int_\mbR \lambda \,dE_\lambda$ if the spectral resolution of $T,$
then the (complex) measure $\mu(a,b) := \tau(A E_{(a,b)} B)$ is countably additive (and has finite variation).
\end{lemma}
\begin{proof} Since the spectral resolution of a self-adjoint operator is strong operator $\sigma$-additive (see e.g. \cite[VIII.3]{RS1}),
the assertion of the lemma follows from Lemma~\ref{L: net VA(alpha) converges to VA}.
\end{proof}

\section{Integration of operator-valued functions}\label{SubS: Def of Oper.V. integral}
\begin{lemma}\label{L: unit ball of E is Polish}
   An invariant operator ideal $\clE$ has property (F) if and only if
   the unit ball of $\clE$ endowed with $so^*$-topology is a complete separable metrisable space.
\end{lemma}
\begin{proof} The "if" part is evident.
Since $\hilb$ is separable,
the unit ball $(\clBH_1, so^*)$ of $\clBH$ is a metrisable space \cite[Proposition I.3.1]{DixvNa}.
Hence, the unit ball $(\clE_1, so^*)$ of $\clE$ is also metrisable.
Since $\hilb$ is separable the unit ball $(\clBH_1, so^*)$ is also separable.
Thus, every subset of $(\clBH_1, so^*)$ is separable \cite[I.6.12]{DS},
and in particular $\clE_1.$ Since the unit ball $(\clBH_1, so^*)$ is complete \cite[Prop. 2.4.1]{BR1},
the property (F) of $\clE$ implies that $(\clE_1, so^*)$ is also complete.
\end{proof}

Let $(S,\Sigma,\nu)$ be a finite measure space and $\clE$ be an invariant operator ideal with property (F).
A bounded function $f\colon(S,\nu) \to \clE$ will be called
\\ \mbox{\quad}  $(i)$ weakly measurable if, for any $\eta,\xi \in \hilb,$
        the function $\scal{f(\cdot)\eta}{\xi}$ is measurable;
\\ \mbox{\quad}  $(ii)$ \ww measurable if, for all $\eta \in \hilb,$
        the functions $f(\cdot)\eta, f(\cdot)^*\eta\colon(S,\nu) \to \hilb$ are Bochner measurable from $S$ into $\hilb;$
\\ \mbox{\quad}  $(iii)$ $so^*$-measurable if there exist a sequence of simple
       (finitely-valued) measurable functions $f_n \colon S \to \clE$ such that $f_n(\sigma) \to f(\sigma)$
       in the $so^*$-topology for a.\,e. $\sigma \in S.$

\begin{prop}\label{P: two def equiv} If $\clE$ has property (F), then,
for any $\clE$-bounded function $f\colon (S,\nu) \to \clE,$ the following conditions are equivalent.
\\ \mbox{\quad}  (i) $f$ is weakly measurable,
\\ \mbox{\quad}  (ii) $f$ is \ww measurable,
\\ \mbox{\quad}  (iii) $f$ is $so^*$-measurable.
\end{prop}
\begin{proof} The implications (iii) $\then$ (ii) $\then$ (i) are evident (and do not depend on property (F)).
That (i) $\then$ (iii) follows from Lemma~\ref{L: unit ball of E is Polish} and \cite[Propositions 1.9 and 1.10]{VTCh}.
\end{proof}

We denote the set of all $\norm{\cdot}$-bounded \ww measurable functions $f\colon S \to \clE$ by $\LsInf(S,\nu;\clE).$
Examples of such functions are bounded $\norm{\cdot}$-Bochner-measurable functions and,
in the case that $S$ is a locally compact space, all $so^*$-continuous bounded functions.

The following lemma is a simple consequence of the previous proposition (cf.\cite[Lemmas 5.5, 5.6]{dPS04FA}).
\begin{lemma}\label{L: the measurables form an algebra} \cite{dPS04FA} (i) The set $\LsInf(S,\nu;\clE)$ is a $*$-algebra; \\
(ii) if
  $\phi \in B_\mbR(\mbR),$ $f \in \LsInf(S,\nu;\clBH_{sa}),$ then $\phi(f) \in \LsInf(S,\nu).$
\end{lemma}

For any bounded function $f \in \LsInf(S,\nu;\clE),$
we define the integral $\int_S f(\sigma)\,d\nu(\sigma)$ by the formula
\begin{gather}\label{F: def of int}
  \brs{\int_S f(\sigma)\,d\nu(\sigma)}\eta = \int_S f(\sigma)\eta\,d\nu(\sigma),
\end{gather}
where the last integral is a Bochner integral.
Evidently, such an integral exists and it is a bounded linear operator
with (uniform) norm less or equal to $\abs{\nu}\norm{f}_\infty.$

\begin{lemma}\label{L: Pointwise linit is meas} If $\clE$ has property (F), and if
  the sequence $f_n \in \LsInf(S,\nu;\clE),$ $n=1,2,\ldots$ is $\clE$-bounded and
  $\nu$-a.\,e. converges to $f\colon S \to \clBH$ in the $so^*$-topology, then
  $f \in \LsInf(S,\nu;\clE).$
\end{lemma}
\begin{proof} We have that, for any $\eta \in \hilb,$ the sequence $f_n(\sigma)\eta$
converges to $f(\sigma)\eta$ for $\nu$-a.e. $\sigma \in S.$
Since the $\hilb$-valued functions $f_n(\cdot)\eta$ are Bochner measurable and since the pointwise
limit of a sequence of Bochner measurable functions is a Bochner measurable function, we have that $f \in \LsInf(S,\nu).$
That $f(\sigma) \in \clE$ for a.\,e. $\sigma \in S$ follows from property (F).
\end{proof}

\begin{lemma}\label{L: int in clE} If $\clE$ has property (F), $f \in \LsInf(S,\nu;\clE)$
and if $f$ is uniformly $\clE$-bounded, then $\int_S f\,d\nu \in \clE.$
\end{lemma}
\begin{proof} By Proposition~\ref{P: two def equiv}, we can choose a sequence
of simple functions $f_n \in \LsInf(S,\nu;\clE)$ converging a.\,e. in $so^*$-topology to $f.$
Evidently, $A_n:=\int_S f_n\,d\nu \in \clE$ for all $n\in\mbN.$ By the definition (\ref{F: def of int})
of operator-valued integral, the sequence $\set{A_n}_{n=1}^\infty$ converges to $\int_S f\,d\nu$
in the $so^*$-topology by the Lebesgue Dominated
Convergence Theorem for the Bochner integral. That $\int_S f\,d\nu \in \clE$ now follows from
the property (F) of $\clE.$
\end{proof}
\begin{cor} Under the assumptions of Lemma~\ref{L: Pointwise linit is meas}, we have
$$
  \int_S f_n \, d\nu \to \int_S f \, d\nu
$$
in the $so^*$-topology.
\end{cor}
\begin{lemma}\label{L: int is oper linear}
  For any $A \in \clBH$ and $B \in \LsInf(S,\nu;\clE)$
    $$A\int_S B(\sigma)\,d\nu(\sigma) = \int_S AB(\sigma)\,d\nu(\sigma)$$
\end{lemma}
The lemma follows directly from \cite[Corollary V.5.2]{Yo}.
\begin{lemma}\label{L: Fubini} If $(S_i,\Sigma_i,\nu_i), \ i=1,2$ are two finite measure spaces and
if $f \in \LsInf(S_1 \mytimes S_2,\nu_1 \mytimes \nu_2),$ then $f(\cdot,t) \in \LsInf(S_1,\nu_1)$ for almost all
$t \in S_2$ and
  \begin{gather}\label{F: FubiniEq}
    \int_{S_2} \int_{S_1} f(s,t) \,d\nu_1(s) \,d\nu_2(t) = \int_{S_1\times S_2} f(s,t) \,d(\nu_1\mytimes\nu_2)(s,t).
  \end{gather}
\end{lemma}
\begin{proof}
Since $f(\cdot,\cdot)$ is integrable, for any $\eta \in \hilb$ there exists a $\nu_2$-measure zero
set $A_\eta \subset S_2$ such that for all $t \notin A_\eta$ the function $f(\cdot,t)\eta$ is Bochner integrable
(see \cite[Theorem III.11.13]{DS}). If $\seque \xi$ is an orthonormal basis in $\hilb$ and $A = \bigcup_{j=1}^\infty A_{\xi_j},$ then
$\nu_2(A) = 0$ and, for any $\eta \in \hilb$ and $t \notin A,$ we have
$$
  f(\cdot,t)\eta = \sums{j=1}^\infty c_n f(\cdot,t)\xi_n,
$$
where $\eta = \sums{j=1}^\infty c_n \xi_n.$ Since linear combinations and uniformly bounded pointwise limits
of sequences
of Bochner integrable functions on the measure space $(S,\nu)$ are Bochner integrable (by the Lebesgue Dominated Convergence
Theorem), it follows that $f(\cdot,t)\eta$ is integrable for $t\notin A.$ Similarly, there exists a $\nu_2$-measure zero set
$A'$ such that $f(\cdot,t)^*\eta$ is integrable for all $\eta \in \hilb$ and $t\notin A'.$
Hence, $f(\cdot,t)$ is integrable for all $t \notin A \cup A'$ and the operator-valued function
$g(t) := \int_{S_1} f(s,t)\,d\nu_1(s)$ is well-defined. Now, the integral $g(t)\eta = \int_{S_1} f(s,t)\eta\,d\nu_1(s)$
exists and is equal to $\int_{S_1\mytimes S_2} f(s,t)\eta\,d(\nu_1\mytimes \nu_2)(s,t)$ by
Fubini's theorem for the Bochner integral of $\hilb$-valued functions \cite[Theorem III.11.13]{DS}. The latter
means that the equality (\ref{F: FubiniEq}) holds.
\end{proof}
\begin{lemma} \label{L: integral for vNa} If $f \in \LsInf(S,\nu; \clN),$ then
\\ (i) $X:=\int_S f(\sigma)\,d\nu(\sigma)$ belongs to $\clN;$
\\ (ii) $X$ as an element of the $W^*$-algebra $\clN$ does not depend on any representation of $\clN.$
\end{lemma}
\begin{proof} (i) Let $A' \in \clN'.$ Then by Lemma~\ref{L: int is oper linear}
\begin{multline*}
  A'X\eta = \int_S A'f(\sigma)\eta\,d\nu(\sigma)
    \\ = \int_S f(\sigma)A'\eta\,d\nu(\sigma) = \int_S f(\sigma)\,d\nu(\sigma)A'\eta = XA'\eta
\end{multline*}
for any $\eta \in \hilb.$ Hence, $X \in \clN.$ \\
(ii) This follows from the fact that two representations of a \vNa\ can be obtained
from each other by ampliation, reduction and spatial isomorphism \cite{DixvNa},
since for each of these isomorphisms the claim is evident.
\end{proof}
\begin{lemma} \label{L: tau int = int tau} If $(\clN,\tau)$ is a semifinite
\vNa, if $f \in \LsInf(S,\nu; \LpN{1})$ and if $f$ is uniformly $\LpN{1}$-bounded,
then $X := \int_S f(\sigma)\,d\nu(\sigma) \in \LpN{1},$ the function $\tau(f(\cdot))$ is measurable and
$$
  \taubrs{\int_S f(\sigma)\,d\nu(\sigma)} =   \int_S \taubrs{f(\sigma)}\,d\nu(\sigma).
$$
\end{lemma}
\begin{proof}
Lemma~\ref{L: int in clE} implies that $X \in \LpN{1},$ so that the left hand side of the equality above makes sense.

By \cite[Lemma 5.9]{dPS04FA}, the function $\tau(f(\cdot))$ is measurable.

By linearity and by Lemma~\ref{L: the measurables form an algebra}(i), we can assume that $f(\cdot) \geq 0.$
By Lemma~\ref{L: integral for vNa}(ii), we can assume
that $\clN$ acts on $L^2(\clN,\tau)$ in the left regular representation. Let $E$ be an arbitrary
$\tau$-finite projection from $\clN.$ Then $E \in L^{2}(\clN,\tau)$ and by the definition
(\ref{F: def of int}) of the operator-valued integral
$$
  XE = \int_S f(\sigma) E\,d\nu(\sigma),
$$
where the right hand side is a Bochner integral in $L^2(\clN,\tau).$
Since $E$ is $\tau$-finite, the convergence in $L^2(\clN,\tau)$ of the Bochner integral
implies convergence in $L^1(\clN,\tau),$ so that we have
$$
  \tau(XE) = \int_S \tau(f(\sigma) E) \,d\nu(\sigma).
$$
Now, normality of the trace $\tau$ and the dominated convergence theorem imply that
$$
  \tau(X) = \int_S \tau(f(\sigma)) \,d\nu(\sigma).
$$
\end{proof}

\section{Multiple operator integrals}\label{S: MOI section}
In this section, we define multiple operator integrals of the form
$$
  \int_{\mbR^{n+1}} \phi(\lambda_0,\ldots,\lambda_n)\,dE^{\HH_0}_{\lambda_0}\VV_1
   \,dE^{\HH_1}_{\lambda_1}\VV_2\,dE^{\HH_2}_{\lambda_2}\ldots \VV_{n}\,dE^{\HH_n}_{\lambda_n}.
$$
We denote by $B(\mbR^{n+1})$ the set of all bounded Borel functions on $\mbR^{n+1}.$
Throughout this section, we consider the set of those functions $\phi \in B(\mbR^{n+1})$
which admit a representation of the form
\begin{gather}\label{F: BS representation}
  \phi(\lambda_0,\lambda_1,\ldots,\lambda_n) = \int_{S} \alpha_0(\lambda_0,\sigma)\ldots\alpha_n(\lambda_n,\sigma)\,d\nu(\sigma),
\end{gather}
where $(S,\nu)$ is a finite measure space and $\alpha_0,\ldots,\alpha_n$ are bounded Borel functions on $\mbR\times S.$
Similar representations (for the case $n=1$) were discussed in \cite{dPS04FA}.

\begin{defn}\label{D: def of MOI}
   For arbitrary self-adjoint operators $\HH_0,\ldots,\HH_n$ on the Hilbert space $\hilb,$ bounded operators $\VV_1,\ldots,\VV_{n}$ on $\hilb$
   and any function $\phi \in B(\mbR^{n+1})$ which admits a representation given by (\ref{F: BS representation}),
   the multiple operator integral $T^{\HH_0,\ldots,\HH_n}_\phi(\VV_1,\ldots,\VV_{n})$ is defined as
\begin{gather}\label{F: Def of MOI}
  T^{\HH_0,\ldots,\HH_n}_\phi(\VV_1,\ldots,\VV_{n}) := \int_{S} \alpha_0(\HH_0,\sigma)\VV_1\ldots \VV_{n}\alpha_n(\HH_n,\sigma)\,d\nu(\sigma),
\end{gather}
where the integral is taken in the sense of definition (\ref{F: def of int}).
\end{defn}
\begin{rems} By \cite[Lemma 5.13]{dPS04FA} and Lemma~\ref{L: the measurables form an algebra}(i)
applied to $\clE = \clBH,$
the function $\sigma \mapsto \alpha_0(\HH_0,\sigma)\VV_1\ldots \VV_{n}\alpha_n(\HH_n,\sigma)$
is \ww measurable and therefore the integral above exists.
\end{rems}
\begin{lemma}\label{L: Def of MOI is correct}
The multiple operator integral in Definition~\ref{D: def of MOI} is well-defined
in the sense that it does not depend on the representation (\ref{F: BS representation}) of $\phi.$
\end{lemma}
\begin{proof} We first prove that if the operators $\VV_1,\ldots,\VV_n$ are all
one-dimensional, then the right hand side of (\ref{F: Def of MOI})
does not depend on the representation of $\phi$ given by (\ref{F: BS representation}).

For $\eta,\xi \in \hilb,$ we denote by $\theta_{\eta,\xi}$ the one-dimensional operator
defined by formula
$\theta_{\eta,\xi}\zeta = \scal{\eta}{\zeta}\xi, \ \zeta \in \hilb.$ It is clear that $\Tr(\theta_{\eta,\xi}) = \scal{\eta}{\xi}, \ $
$A\theta_{\eta,\xi} = \theta_{\eta,A\xi}$ for any $A \in \clBH$ and that
$\theta_{\eta_1,\xi_1}\ldots \theta_{\eta_n,\xi_n} = \scal{\eta_1}{\xi_2}\ldots\scal{\eta_{n-1}}{\xi_n}\theta_{\eta_n,\xi_1}.$

Let $\VV_j = \theta_{\eta_j,\xi_j}, \ j = 0,\ldots, n.$ Then
\begin{align*}
   E & :=  \Tr\brs{\VV_0 \int_{S} \alpha_0(\HH_0,\sigma)\VV_1\ldots \VV_{n}\alpha_n(\HH_n,\sigma)\,d\nu(\sigma)}
  \\ &= \Tr\int_{S} \VV_0\alpha_0(\HH_0,\sigma)\VV_1\ldots \VV_{n}\alpha_n(\HH_n,\sigma)\,d\nu(\sigma)
  \\ &= \Tr\int_{S} \theta_{\eta_0,\xi_0}\alpha_0(\HH_0,\sigma)\theta_{\eta_1,\xi_1}\ldots \theta_{\eta_n,\xi_n}\alpha_n(\HH_n,\sigma)\,d\nu(\sigma)
  \\ &= \int_{S} \Tr\brs{\theta_{\eta_0,\xi_0}\alpha_0(\HH_0,\sigma)\theta_{\eta_1,\xi_1}\ldots \theta_{\eta_n,\xi_n}\alpha_n(\HH_n,\sigma)}\,d\nu(\sigma)
  \\ &= \int_{S} \Tr\brs{\alpha_0(\HH_0,\sigma)\theta_{\eta_1,\xi_1}\ldots \theta_{\eta_n,\xi_n}\alpha_n(\HH_n,\sigma)\theta_{\eta_0,\xi_0}}\,d\nu(\sigma)
  \\ &= \int_{S} \Tr\brs{\theta_{\eta_1,\alpha_0(\HH_0,\sigma)\xi_1}\ldots \theta_{\eta_n,\alpha_{n-1}(\HH_{n-1},\sigma)\xi_n}\theta_{\eta_0,\alpha_n(\HH_n,\sigma)\xi_0}}\,d\nu(\sigma)
  \\ &= \int_{S} \scal{\eta_0}{\alpha_0(\HH_0,\sigma)\xi_1}\scal{\eta_1}{\alpha_1(\HH_1,\sigma)\xi_2}\ldots \scal{\eta_n}{\alpha_{n}(\HH_{n},\sigma)\xi_0}\,d\nu(\sigma).
\end{align*}
Now, since $\scal{\eta}{\alpha(\HH)\xi} = \int_\mbR \alpha(\lambda) \scal{\eta}{dE^\HH_\lambda\xi},$
we have that
$$
  E=\int_{S} \int_\mbR \alpha_0(\lambda_0,\sigma) \scal{\eta_0}{dE^{\HH_0}_{\lambda_0}\xi_1}
    \ldots \int_\mbR \alpha_n(\lambda_n,\sigma) \scal{\eta_n}{dE^{\HH_n}_{\lambda_n}\xi_0} \,d\nu(\sigma).
$$
Since the measure $\scal{\eta}{dE_\lambda\xi}$ has finite total variation, Fubini's theorem implies
\begin{align*}
  E& = \int_{S} \brs{\int_{\mbR^{n+1}} \alpha_0(\lambda_0,\sigma)\ldots\alpha_n(\lambda_n,\sigma)
  \scal{\eta_0}{dE^{\HH_0}_{\lambda_0}\xi_1} \ldots \scal{\eta_n}{dE^{\HH_n}_{\lambda_n}\xi_0}}\,d\nu(\sigma)
  \\ & = \int_{\mbR^{n+1}} \brs{\int_{S} \alpha_0(\lambda_0,\sigma)\ldots\alpha_n(\lambda_n,\sigma)\,d\nu(\sigma)}
  \scal{\eta_0}{dE^{\HH_0}_{\lambda_0}\xi_1}  \ldots \scal{\eta_n}{dE^{\HH_n}_{\lambda_n}\xi_0}
  \\ & = \int_{\mbR^{n+1}} \phi(\lambda_0,\ldots,\lambda_n)
     \scal{\eta_0}{dE^{\HH_0}_{\lambda_0}\xi_1} \ldots \scal{\eta_n}{dE^{\HH_n}_{\lambda_n}\xi_0}.
\end{align*}
We recall that, if $A, B$ are bounded operators, then $A=B$ if and only if the equality $\Tr(\VV A) = \Tr(\VV B)$
holds for all one-dimensional operators $\VV.$
It now follows immediately that the multiple operator integral does not
depend on the representation (\ref{F: BS representation}) of $\phi$ in the case that the operators $\VV_1,\ldots,\VV_n$ are one-dimensional.

By linearity, it follows that the definition of multiple operator integral
does not depend on the representation (\ref{F: BS representation}) in the case of
finite-dimensional operators $\VV_1,\ldots,\VV_n.$
Since every bounded operator is an $so$-limit of a sequence of finite-dimensional operators,
the claim follows from Proposition~\ref{P: Continuity of MOI}.
\end{proof}
\begin{lemma}\label{L: vNa and MOI}
   If $\clN$ is a \vNa, if $\HH_0,\ldots,\HH_n$ are self-adjoint operators affiliated with $\clN$
   and if $\VV_1,\ldots,\VV_n \in \clN,$ then $T^{\HH_0,\ldots,\HH_n}_\phi(\VV_1,\ldots,\VV_n) \in \clN.$
\end{lemma}
This follows from Lemma~\ref{L: integral for vNa}.

The following observation is a direct consequence of Lemma
\ref{L: integral formula for n-th divided difference} and Definition~\ref{D: def of MOI}.
\begin{lemma} \label{L: MOI for f[n]}
 If $f \in \CPlus{n},$ then
  \begin{multline}
    T^{\HH_0,\ldots,\HH_n}_{f^{[n]}}(\VV_1,\ldots,\VV_n) = \\
     \int_{\Pin{n}} e^{i(s_0-s_1)\HH_0} \VV_1 e^{i(s_{1}-s_2)\HH_1}\VV_2
     \ldots \VV_n e^{is_nB_n} \,d\nunf{n}(s_0,\ldots,s_n).
  \end{multline}
\end{lemma}
\begin{lemma}\label{L: MOI clE is subset of clE} If $\clE$ is an invariant operator ideal with property (F)
  and if one of the operators $\VV_1,\ldots,\VV_n$ belongs to $\clE,$ then
  $$
    T^{\HH_0,\ldots,\HH_n}_{\phi}(\VV_1,\ldots,\VV_n) \in \clE.
  $$
  In case that $n=2,$ this yields
  $$
    \norm{T^{\HH_1,\HH_2}_{\phi}}_{\clE \to \clE} \leq \norm{\phi},
  $$
  where (see \cite{dPS04FA})
  \begin{multline}
    \norm{\phi} = \inf \Big\{\int_S \norm{\alpha(\cdot,\sigma)}_\infty\norm{\beta(\cdot,\sigma)}_\infty\,d\nu(\sigma)
      \colon \phi(\lambda,\mu)
      \\ = \int_S \alpha(\lambda,\sigma)\beta(\mu,\sigma)\,d\nu(\sigma)\Big\}.
  \end{multline}
\end{lemma}
\begin{proof}
    Follows from Lemmas~\ref{L: the measurables form an algebra}(i) and~\ref{L: int in clE}.
\end{proof}
\begin{rems}
  If $\VV \in \LpN{2}$ and if $n=2,$ then the preceding definition
  coincides with the definition of double operator integral as
  a spectral integral given in \cite{BS73DOI} and \cite{dPS04FA}.
\end{rems}
\begin{cor}\label{C: trace of MOI}
 If $\VV_1,\ldots,\VV_n \in \clN,$ $\VV_j \in \LpN{1}$ for some $j = 1,\ldots,n,$ $\HH_0, \ldots, \HH_n$
 are self-adjoint operators affiliated with $\clN,$ $\phi \in B(\mbR^{n+1})$ and $\phi(\lambda_0,\ldots,\lambda_n)$
 admits the representation (\ref{F: BS representation}), then
  \begin{multline}
    \taubrs{T^{\HH_0,\ldots,\HH_n}_\phi(\VV_1,\ldots,\VV_n)}
    \\ = \int_S \taubrs{\alpha_0(\HH_0,\sigma)\VV_1\alpha_1(\HH_1,\sigma)\ldots \VV_n\alpha_n(\HH_n,\sigma)}\,d\nu(\sigma)
  \end{multline}
\end{cor}
\begin{proof} It is enough to note that the operator-valued function
 $$\sigma \mapsto \alpha_0(\HH_0,\sigma)\VV_1\alpha_1(\HH_1,\sigma)\ldots \VV_n\alpha_n(\HH_n,\sigma)$$
is \ww measurable by \cite[Lemma 5.11]{dPS04FA} and Lemma~\ref{L: the measurables form an algebra}(i), so that we can apply
Lemma~\ref{L: tau int = int tau}.
\end{proof}

\begin{prop}\label{P: Continuity of MOI} (i)
If a sequence of self-adjoint operators $\VV^{(k_j)}_j \in \clBH, j = 1,\ldots,n,$
converges to $\VV_j \in \clBH$ in the $so$-topology (respectively, norm topology) as $k_j \to \infty,$ then
$$
  T^{\HH_0,\ldots,\HH_n}_\phi (\VV^{(k_1)}_1,\ldots,\VV^{(k_n)}_n) \to T^{\HH_0,\ldots,\HH_n}_\phi (\VV_1,\ldots,\VV_n)
$$
in the $so$-topology (respectively, norm topology) as $k_1, \ldots, k_n \to \infty.$

(ii)
If a sequence of self-adjoint operators $\HH^{(k_j)}_j, j = 0,\ldots,n$
resolvent strongly converges to $\HH_j$ as $k_j \to \infty$ and $\VV_1,\ldots,\VV_n \in \clBH,$ then
$$
  T^{\HH_0^{(k_0)},\ldots,\HH_n^{(k_n)}}_\phi (\VV_1,\ldots,\VV_n) \to T^{\HH_0,\ldots,\HH_n}_\phi (\VV_1,\ldots,\VV_n)
$$
in the $so$-topology as $k_0, \ldots, k_n \to \infty.$
\end{prop}
\begin{proof} We prove the part (ii), the proof of part (i) is similar (and simpler). Suppose that
$$
  \phi(\lambda_0,\ldots,\lambda_n) = \int_S \alpha_0(\lambda_0,\sigma)\ldots\alpha_n(\lambda_n,\sigma)\,d\nu(\sigma)
$$
is a representation of $\phi$ given by (\ref{F: BS representation}).
Since $\alpha(\cdot,\sigma)$ is a bounded function for every $\sigma \in S,$ the operators $\alpha(\HH_j^{(k_j)},\sigma)$ converge to
$\alpha(\HH_j,\sigma)$ in the $so$-topology \cite[Theorem VIII.20(b)]{RS1}.
Since multiplication of operators is jointly continuous in the $so$-topology
on the unit ball of $\clN$ \cite[Proposition 2.4.1]{BR1},
the operator $\alpha(\HH_0^{(k_0)},\sigma)\VV_1 \ldots \VV_n \alpha(\HH_n^{(k_n)},\sigma)$ converges in the $so$-topology to
$\alpha(\HH_0,\sigma)\VV_1 \ldots \VV_n \alpha(\HH_n,\sigma), \ \sigma \in S.$
Now, an application of the Dominated Convergence Theorem for the Bochner integral
of $\hilb$-valued functions \cite[Corollary III.6.16]{DS}
completes the proof.
\end{proof}
This new definition of multiple operator integral enables us to give
a simple proof of the following
\begin{prop}\label{P: Prop of MOI} The multiple operator integral has the properties: \\
(i) if $\phi_1$ and $\phi_2$ admit a representation of the type given in (\ref{F: BS representation}),
then so does $\phi_1+\phi_2$ and
   \begin{gather}\label{F: MOI of phi 1 plus phi 2}
     T^{\HH_1,\ldots,\HH_n}_{\phi_1+\phi_2} = T^{\HH_1,\ldots,\HH_n}_{\phi_1} + T^{\HH_1,\ldots,\HH_n}_{\phi_2};
   \end{gather}
(ii) in the case of double operator integrals, if $\phi_1$ and $\phi_2$ admit a representation
of the type given in (\ref{F: BS representation}), then so does $\phi_1\phi_2$ and
   $$
     T^{\HH_1,\HH_2}_{\phi_1\phi_2} = T^{\HH_1,\HH_2}_{\phi_1}T^{\HH_1,\HH_2}_{\phi_2}.
   $$
\end{prop}
\begin{proof}
   (i) If we take representations of the form (\ref{F: BS representation}) with $(S_1,\nu_1)$ and $(S_2,\nu_2)$ for $\phi_1$ and $\phi_2$
   and put $(S,\nu) = (S_1,\nu_1) \sqcup (S_2,\nu_2)$ for $\phi_1+\phi_2$ with evident definition of $\alpha_1, \alpha_2, \ldots,$
   then the equality (\ref{F: MOI of phi 1 plus phi 2}) follows from Definition~\ref{D: def of MOI}. Here $\sqcup$
   denotes the disjoint sum of measure spaces.

   (ii) If $$\phi_j(\lambda_1,\lambda_2) = \int_{S_1} \alpha_j(\lambda_1,\sigma_1)\beta_j(\lambda_2,\sigma_1)\,d\nu_j(\sigma_1), \ j = 1,2,$$
   set
   $$\phi(\lambda_1,\lambda_2) = \int_{S} \alpha(\lambda_1,\sigma)\beta(\lambda_2,\sigma)\,d\nu(\sigma),$$
   where
   $$
     (S,\nu) = (S_1,\nu_1) \mytimes (S_2,\nu_2)
   $$ and
   $$
     \quad \alpha(\lambda,\sigma) = \alpha_1(\lambda,\sigma_1)\alpha_2(\lambda,\sigma_2), \quad
    \beta(\lambda,\sigma) = \beta_1(\lambda,\sigma_1)\beta_2(\lambda,\sigma_2).
   $$
     Consequently,
   \begin{gather*}
     T^{\HH_1,\HH_2}_{\phi_1}\brs{T^{\HH_1,\HH_2}_{\phi_2}(\VV)} = \int_{S_1} \alpha_1(\HH_1,\sigma_1)T^{\HH_1,\HH_2}_{\phi_2}(\VV)\beta_1(\HH_2,\sigma_1)\,d\nu_1(\sigma_1)
     \\ = \int_{S_1} \alpha_1(\HH_1,\sigma_1)\brs{\int_{S_2} \alpha_2(\HH_1,\sigma_2)\VV\beta_2(\HH_2,\sigma_2)\,d\nu_2(\sigma_2)}\beta_1(\HH_2,\sigma_1)\,d\nu_1(\sigma_1).
   \end{gather*}
   Now, Lemma~\ref{L: int is oper linear} and Fubini's theorem (Lemma~\ref{L: Fubini}) imply
   \begin{multline*}
     T^{\HH_1,\HH_2}_{\phi_1}\brs{T^{\HH_1,\HH_2}_{\phi_2}(\VV)}
       = \int_{S_1\mytimes S_2} \alpha_1(\HH_1,\sigma_1)\alpha_2(\HH_1,\sigma_2)\VV
          \\ \times \beta_2(\HH_2,\sigma_2)\beta_1(\HH_2,\sigma_1)\,d(\nu_1\mytimes\nu_2)
     (\sigma_1,\sigma_2) = T^{\HH_1,\HH_2}_{\phi_1\phi_2}(\VV).
   \end{multline*}
\end{proof}
\section{Higher order Fr\'echet differentiability}\label{S: Dal-Krein formula}
We note that, by Stone's theorem \cite[Theorem VIII.7]{RS1}
and joint continuity of multiplication of operators (from the unit ball)
in the $so$-topology \cite[Proposition 2.4.1]{BR1}
all operator-valued integrals occurring in this and subsequent sections
are defined as in section~\ref{SubS: Def of Oper.V. integral}.
\begin{lemma}\label{L: Operator Fourier transform}
If $A$ is a self-adjoint (possibly unbounded) operator
on a~Hilbert space~$\hilb$  and if $f$ is a function on $\mbR$ such that $f \in \CPlus 1,$
then
$$
  f(A) = (2\pi)^{-1/2} \int_\mbR e^{isA}\,m_f(ds).
$$
\end{lemma}
The proof is a simple application of Fubini's theorem. See \cite[Theorem 3.2.32]{BR1}
\begin{lemma}\label{L: Duhamel's formula} (Duhamel's formula).
  If $B$ is an unbounded self-adjoint operator on a~Hilbert space~$\hilb,$
  if $\VV$ is a bounded self-adjoint operator on~$\hilb$ and if $A=B+\VV,$ then
  \begin{gather}\label{F: Duhamel formula II}
    e^{isA} - e^{is\HH} = \int_0^s e^{i(s-t)A} i(A-B) e^{it\HH}\,dt.
  \end{gather}
\end{lemma}
 \begin{proof}
   Let $F(t) = e^{itA}e^{-it\HH}.$ Taking derivative of $F(t)$ in the $so$-topology gives
   $$
     F'(t) = iA e^{itA}e^{-it\HH}+ e^{itA}(-iB)e^{-it\HH} = e^{itA} i(A-B) e^{-it\HH}.
   $$
   So,
   $$
     \int_0^s e^{itA} i(A-B) e^{-it\HH}\,dt = F(s) - F(0) = e^{isA}e^{-is\HH} - 1.
   $$
   Multiplying the last equality by $e^{is\HH}$ from the right gives (\ref{F: Duhamel formula II}).
 \end{proof}
\begin{thm}\label{T: DalKreinIIa}
  Let $\clN$ be a \vNa. Suppose that $\HH=\HH^*$ is affiliated with $\clN,$ that
  $\VV \in \clN$ is self-adjoint and set $A=\HH+\VV.$ If $f \in \CPlus 1,$ then
    \begin{gather*}
      f(A) - f(\HH) = T^{A,\HH}_{\psif{f}}(\VV).
    \end{gather*}
\end{thm}
\begin{proof}
It follows from Lemma~\ref{L: Operator Fourier transform} that
\begin{gather*}
  f(A) - f(\HH) = \frac {1}{\sqrt{2\pi}} \int_\mbR (e^{isA}-e^{is\HH})\,m_f(ds).
\end{gather*}
Hence, by Lemma~\ref{L: Duhamel's formula},
\begin{gather}\label{F: first int}
  f(A) - f(\HH) = \frac {i}{\sqrt{2\pi}} \int_\mbR m_f(ds) \int_0^s e^{i(s-t)A} \VV e^{it\HH}\,dt.
\end{gather}
Since $f\in \CPlus 1,$ by Lemma~\ref{L: (Pi,nu) is a finite measure space} and Fubini's theorem (Lemma~\ref{L: Fubini}),
the repeated integral can be replaced by a double integral, so that
\begin{multline}\label{F: sec int}
  f(A) - f(\HH) = \frac {i}{\sqrt{2\pi}} \intint_{\Pi} e^{i(s-t)A} \VV e^{it\HH}\,m_f(ds) \,dt
             \\ = \intint_{\Pi} e^{i(s-t)A} \VV e^{it\HH}\,d\nu_f(\sigma).
\end{multline}
It now follows from Lemma~\ref{L: MOI for f[n]}
that $f(A) - f(\HH) = T^{A,\HH}_{\psif{f}}(\VV).$
\end{proof}
\begin{rems} The preceding formula is due to Birman-Solomyak \cite{BS73DOI}.
It is similar to \cite[Corollary 7.2]{dPSW02FA}, which applies to a wider class of functions
but is restricted to bounded operators in a semifinite \vNa\ $\clN.$
\end{rems}
Let $X$ be a topological vector space, $\clE$ be a normed space embedded in $X.$
Let $x \in X$ and $f \colon x+\clE \to f(x) + \clE.$ The function $f$ is called
affinely Fr\'echet differentiable at $x$ along $\clE$ if there exists
a (necessarily unique) bounded operator $L \colon \clE \to \clE$ such that
$$
  f(x+h) - f(x) = L(h) + r(x,h),
$$
where $\norm{r(x,h)}_{\clE} = o(\norm{h}_{\clE}).$ We write $L = D_\clE f(x).$
\begin{thm}\label{T: DalKreinIIb}
  Let $\clN$ be a von Neumann algebra, acting in a Hilbert space $\hilb$.
  Let $\HH=\HH^*$ be affiliated with $\clN$ and let $\VV \in \clE_{sa},$ where $\clE$
  is an invariant operator ideal over $\clN$ with property (F).
  If $f \in \CPlus{2},$ then the function $f\colon \HH' \in \HH+\clE_{sa} \mapsto f(\HH') \in f(\HH)+\clE_{sa}$
  is affinely Fr\'echet differentiable along $\clE_{sa}$ and $D_{\clE} f(\HH) = T^{\HH,\HH}_{\psif{f}}.$
  The function $X \mapsto D_\clE f(\HH+X)$ is continuous in the norm of $\clE$ and satisfies the estimate
    \begin{multline} \label{F: norm estimate}
      \norm{D_\clE f(\HH+X)(\VV)-D_\clE f(\HH)(\VV)}_{\clE}
      \\ \leq \norm{m_{f''}} \norm{\VV}_\clE \norm{X}_\clE,
           \ X, \VV \in \clE.
    \end{multline}
\end{thm}
\begin{proof}
By (\ref{F: sec int}) we have, following \cite{Wi},
\begin{multline}\label{F: sec int 2}
  f(\HH+\VV) - f(\HH) = \intint_{\Pi} e^{i(s-t)(\HH+\VV)} \VV e^{it\HH}\,d\nu_f(s,t) \\
            = \intint_{\Pi} e^{i(s-t)\HH} \VV e^{it\HH}\,d\nu_f(s,t)
                + \intint_{\Pi} \brackets{e^{i(s-t)(\HH+\VV)}-e^{i(s-t)\HH}} \VV e^{it\HH}\,d\nu_f(s,t)
              \\  = (I) + (II).
\end{multline}
$(I)$ is equal to $T^{\HH,\HH}_{\psif{f}}(\VV)$ and represents a continuous
linear operator on~$\clE$ (see Lemmas~\ref{L: MOI for f[n]} and~\ref{L: MOI clE is subset of clE}),
so that it will be a Fr\'echet derivative of $f\colon \HH+\clE \to f(\HH)+\clE$
provided it is shown that the second term is $o(\norm{\VV}_\clE).$ Applying
Duhamel's formula (\ref{F: Duhamel formula II}) yields
\begin{gather}\label{F: II}
  (II) = \intint_{\Pi} \brackets{\int_0^{s-t} e^{i(s-t-u)(\HH+\VV)} i \VV e^{iu\HH}\,du} \VV e^{it\HH}\,d\nu_f(s,t).
\end{gather}
Since $f\in \CPlus{2},$ Lemmas~\ref{L: calculus exercise},~\ref{L: Def of MOI is correct} and~\ref{L: MOI for f[n]}
enable us to rewrite (\ref{F: II}) as
\begin{gather*}
  (II) = \intint\!\!\!\!\int_{\Pin{2}} e^{i(s-t)(\HH+\VV)} \VV e^{i(t-u)\HH} \VV e^{iu\HH}\,d\nunf{2}(s,t,u),
\end{gather*}
where $(\Pin{2},\nunf{2})$ is the finite measure space defined in Lemma~\ref{L: (Pi,nu) is a finite measure space}.
The $\clE$-norm of the last expression is estimated by $\abs{\nunf{2}} \norm{\VV} \norm{\VV}_\clE \leq \abs{\nunf{2}}\norm{\VV}^2_\clE.$
So, the function $f\colon \HH+\clE \to f(\HH)+\clE$ is Fr\'echet differentiable and $D_{\clE} f(\HH) = T^{\HH,\HH}_{\psif{f}}.$

The norm continuity of this derivative and the estimate (\ref{F: norm estimate})
follow by a similar argument using Duhamel's formula (\ref{F: Duhamel formula II}).
\end{proof}
\begin{rems}
It follows, in particular, from the preceding theorem via Lemma~\ref{L: MOI clE is subset of clE}
that the operator $\left.T^{\HH,\HH}_{\psif{f}}\right|_\clE$ is a bounded
linear operator on $\clE.$
\end{rems}
\begin{thm}\label{T: n-th Frechet derivative}
Let $\clN$ be a \vNa\ on a Hilbert space $\hilb,$ let $\HH=\HH^*$ be affiliated with $\clN$
and let $\VV_1,\ldots,\VV_n \in \clE_{sa}.$
If $f\in \CPlus{n+1},$ then the function $f\colon \HH' \in \HH+\clE_{sa} \mapsto f(\HH') \in f(\HH)+\clE_{sa}$
is $n$-times affinely Fr\'echet differentiable along $\clE_{sa}$ and
\begin{gather}\label{F: n-th Frechet der}
  D_\clE^n f(\HH)(\VV_1,\ldots,\VV_n) = \sums{\sigma \in P_n} T^{\HH,\ldots,\HH}_{f^{[n]}}(\VV_{\sigma(1)},\ldots,\VV_{\sigma(n)}) \in \clE,
\end{gather}
where $P_n$ is the standard permutation group.
\end{thm}
\begin{proof}
If $n=1$ then this theorem is exactly Theorem~\ref{T: DalKreinIIb}.
Set $\tHH = \HH + \VV_{n+1}.$
By induction we have
\begin{gather*}
  D^n f(\tHH; \VV_1,\ldots, \VV_n) - D^n f(\HH; \VV_1,\ldots, \VV_n) \\
  = \sums{\sigma \in P_n}  \brs{ T_{f^{[n]}}^{\tHH,\tHH,\ldots,\tHH}(\VV_{\sigma(1)},\ldots,\VV_{\sigma(n)})
    - T_{f^{[n]}}^{\HH,\HH,\ldots,\HH}(\VV_{\sigma(1)},\ldots,\VV_{\sigma(n)}) }.
\end{gather*}
A single term of this sum is
\begin{gather*}
  T_{f^{[n]}}^{\tHH,\tHH,\ldots,\tHH}(\VV_{\sigma(1)},\ldots,\VV_{\sigma(n)}) - T_{f^{[n]}}^{\HH,\HH,\ldots,\HH}(\VV_{\sigma(1)},\ldots,\VV_{\sigma(n)})
    \\ = \sums{j=0}^{n} \brs{T_{f^{[n]}}^{\tHH,\ldots,\stackrel{(j)}{\tHH},\HH,\ldots,\HH}(\VV_{\sigma(1)},\ldots,\VV_{\sigma(n)})
    - T_{f^{[n]}}^{\tHH,\ldots,\tHH,\stackrel{(j)}{\HH},\ldots,\HH}(\VV_{\sigma(1)},\ldots,\VV_{\sigma(n)})}.
\end{gather*}
Now, the $j$-th summand is (Lemma~\ref{L: MOI for f[n]})
\begin{gather*}
  T_{f^{[n]}}^{\tHH,\ldots,\stackrel{(j)}{\tHH},\HH,\ldots,\HH}(\VV_{\sigma(1)},\ldots,\VV_{\sigma(n)})
    - T_{f^{[n]}}^{\tHH,\ldots,\tHH,\stackrel{(j)}{\HH},\ldots,\HH}(\VV_{\sigma(1)},\ldots,\VV_{\sigma(n)})
  \\=\int_{\Pin{n}} e^{i(s_0-s_1)\tHH}\VV_{\sigma(1)} \ldots \VV_{\sigma(j)} e^{i(s_{j}-s_{j+1})\tHH}\VV_{\sigma(j+1)} e^{i(s_{j+1}-s_{j+2})\HH} \VV_{\sigma(j+2)}
        \\ \ldots \VV_{\sigma(n)} e^{is_nB}\,d\nunf{n}(s_0,\ldots,s_n)
    \\-\int_{\Pin{n}} e^{i(s_0-s_1)\tHH}\VV_{\sigma(1)} \ldots \VV_{\sigma(j-1)} e^{i(s_{j-1}-s_{j})\tHH}\VV_{\sigma(j)} e^{i(s_{j}-s_{j+1})\HH} \VV_{\sigma(j+1)}
        \\ \ldots \VV_{\sigma(n)} e^{is_nB}\,d\nunf{n}(s_0,\ldots,s_n)
  \\=\int_{\Pin{n}} e^{i(s_0-s_1)\tHH}\VV_{\sigma(1)} \ldots \VV_{\sigma(j)} \brs{e^{i(s_{j}-s_{j+1})\tHH}-e^{i(s_{j}-s_{j+1})\HH}} \VV_{\sigma(j+1)}
        \\ \ldots \VV_{\sigma(n)} e^{is_nB}\,d\nunf{n}(s_0,\ldots,s_n).
\end{gather*}
By Duhamel's formula (Lemma~\ref{L: Duhamel's formula}), we have
\begin{gather*}
  T_{f^{[n]}}^{\tHH,\ldots,\stackrel{(j)}{\tHH},\HH,\ldots,\HH}(\VV_{\sigma(1)},\ldots,\VV_{\sigma(n)})
    - T_{f^{[n]}}^{\tHH,\ldots,\tHH,\stackrel{(j)}{\HH},\ldots,\HH}(\VV_{\sigma(1)},\ldots,\VV_{\sigma(n)})
  \\ = \int_{\Pin{n}} e^{i(s_0-s_1)\tHH}\VV_{\sigma(1)} \ldots \VV_{\sigma(j)}
           \brs{\int_0^{s_{j}-s_{j+1}} e^{iu\tHH} i \VV_{n+1} e^{i(s_{j}-s_{j+1}-u)\HH}\,du}
           \\ \VV_{\sigma(j+1)} \ldots \VV_{\sigma(n)}e^{is_nB}\,d\nunf{n}(s_0,\ldots,s_n).
\end{gather*}
Applying Fubini's theorem (Lemma~\ref{L: Fubini}) we get
\begin{gather}\label{F: integral formula for the difference T f[n](tB...tB...B) - T f[n](tB...B...B)}
  T_{f^{[n]}}^{\tHH,\ldots,\stackrel{(j)}{\tHH},\HH,\ldots,\HH} - T_{f^{[n]}}^{\tHH,\ldots,\tHH,\stackrel{(j)}{\HH},\ldots,\HH}
   \\ = i\int_{\Pin{n}}\int_0^{s_{j}-s_{j+1}}  e^{i(s_0-s_1)\tHH}\VV_{\sigma(1)} \ldots \VV_{\sigma(j)}
        e^{iu\tHH} \VV_{n+1} e^{i(s_{j}-s_{j+1}-u)\HH} \VV_{\sigma(j+1)} \notag
      \\ \ldots \VV_{\sigma(n)} e^{is_nB}\,du\,d\nunf{n}(s_0,\ldots,s_n). \notag
\end{gather}
Hence, it follows from formula (\ref{F: integral formula for the difference T f[n](tB...tB...B) - T f[n](tB...B...B)}),
Lemma~\ref{L: calculus exercise} and the fact that multiple operator integral is well-defined (Lemma~\ref{L: Def of MOI is correct}) that
\begin{gather*}
  T_{f^{[n]}}^{\tHH,\ldots,\stackrel{(j)}{\tHH},\HH,\ldots,\stackrel{(n)}{\HH}}(\VV_{\sigma(1)},\ldots,\VV_{\sigma(n)})
      - T_{f^{[n]}}^{\tHH,\ldots,\tHH,\stackrel{(j)}{\HH},\ldots,\stackrel{(n)}{\HH}}(\VV_{\sigma(1)},\ldots,\VV_{\sigma(n)})
    \\ = T_{f^{[n+1]}}^{\tHH,\ldots,\stackrel{(j)}{\tHH},\HH,\ldots,\stackrel{(n+1)}{\HH}}
        (\VV_{\sigma(1)},\ldots,\VV_{\sigma(j)},\VV_{n+1},\VV_{\sigma(j+1)},\ldots,\VV_{\sigma(n)}).
\end{gather*}
Since the multiple operator integral on the right hand side minus the same
multiple operator integral with the last $\tHH$ replaced by $\HH$
has the order of $o((\max\norm{\VV_j})^{n+2})$ by Duhamel's formula, we see that the theorem is proved.

That the value of the derivative (\ref{F: n-th Frechet der}) belongs to $\clE$
follows from Lemma~\ref{L: MOI clE is subset of clE}.
\end{proof}
The argument of the last proof and Lemma~\ref{L: MOI for f[n]} implies
\begin{cor} Let $\clN$ be a \vNa\ on a Hilbert space $\hilb.$
If $\HH=\HH^*$ is affiliated with $\clN,$ if $\VV \in \clE_{sa}$ and if $f\in \CPlus{n+1},$ then
\begin{multline}
  f(\HH+\VV) - f(\HH) = T^{\HH,\HH}_{f^{[1]}}(\VV) + T^{\HH,\HH,\HH}_{f^{[2]}}(\VV,\VV) +
  \\ \ldots + T^{\HH,\ldots,\HH}_{f^{[n]}}(\VV,\ldots,\VV) + O(\norm{\VV}_\clE^{n+1}).
\end{multline}
\end{cor}
\begin{proof} This corollary is a consequence
  of Theorem~\ref{T: n-th Frechet derivative}
  and Taylor's formula \cite[Theorem 1.43]{Sch}.
\end{proof}
\section{Spectral shift and spectral averaging}\label{S: BirSol in vNa}
The aim of this section is to prove a semifinite extension of
a formula for spectral averaging due to Birman-Solomyak \cite{BS72SM}.

We first recall the following extension of the spectral shift formula of M.\,G.\,Krein from \cite[Theorem 3.1]{ADS}.
\begin{thm}\label{T: Thm 3.1 ADS} If $\HH=\HH^*$ is affiliated with $\clN$ and $\VV=\VV^* \in \LpN{1},$
  then there exists a unique function $\xi=\xi_{\HH+\VV,\HH}(\cdot) \in L^1(\mbR)$
  such that
  \begin{gather*}
    \norm{\xi}_1 \leq \norm{\VV}_1, \quad \int_{-\infty}^\infty \xi(\lambda)\,d\lambda = \tau(\VV), \\
    -\tau(\mathrm{supp}(\VV_-)) \leq \xi(\lambda) \leq \tau(\mathrm{supp}(\VV_+)) \quad \text{for a.e.} \ \lambda \in \mbR
  \end{gather*}
  and for any function $f \in C^1(\mbR),$ whose derivative $f'$ admits the representation
  $$
    f'(\lambda) = \int_{-\infty}^\infty e^{-i\lambda t}\,dm(t), \ \lambda \in \mbR
  $$
  for some finite (complex) Borel measure on $\mbR,$ the operator $f(\HH+\VV) - f(\HH)$
  is $\tau$-trace class and
  $$
    \taubrs{f(\HH + \VV) - f(\HH)} = \int_{-\infty}^\infty f'(\lambda) \xi(\lambda)\,d\lambda.
  $$
\end{thm}
The function $\xi_{\HH+\VV,\HH}(\cdot)$ is called the \emph{Krein spectral shift function} for the pair $(\HH+\VV,\HH).$
\begin{lemma}\label{L: trace of VE(l,s) is measurable}
  If $(\clN,\tau)$ is a semifinite \vNa, if $\HH = \HH^*$ if affiliated with $\clN$
  and $\VV \in \LpN{1},$ then the function $\gamma(\lambda,\rr) = \taubrs{\VV E^{\HH_\rr}_\lambda}$
  is measurable, where $\HH_\rr := \HH+\rr\VV, \ \rr \in [0,1].$
\end{lemma}
\begin{proof} Let $\phi_{\lambda,\eps}$ be a smooth approximation of $\chi_{(-\infty,\lambda]}.$
We note that $\phi_{\lambda,\eps}(\HH) = \phi_{0,\eps}(\HH -\lambda),$ and that the unbounded-operator valued
function $(\lambda,\rr)\in \mbR^2 \mapsto \HH_\rr-\lambda$ is resolvent uniformly continuous \cite[VIII.7]{RS1}.
It follows from \cite[Theorem VIII.20(b)]{RS1} that the function $(\lambda,\rr) \mapsto \phi_{\lambda,\eps}(\HH_\rr)$
is $so$-continuous, so that Lemma~\ref{L: net VA(alpha) converges to VA} implies
that the function $(\lambda,\rr) \mapsto \taubrs{\VV\phi_{\lambda,\eps}(\HH_\rr)}$ is continuous.
Now, since $\phi_{\lambda,\eps} \to \chi_{(-\infty,\lambda]}$ pointwise as $\eps \to 0,$ the operator
$\phi_{\lambda,\eps}(\HH_\rr) \to \chi_{(-\infty,\lambda]}(\HH_\rr)$ in $so$-topology.
Hence, again by Lemma~\ref{L: net VA(alpha) converges to VA},
the function $\taubrs{\VV\chi_{(-\infty,\lambda]}(\HH_\rr)}$ is measurable.
\end{proof}
\begin{thm}\label{T: BirmanSolII} Let $(\clN,\tau)$ be a semifinite von Neumann algebra on a Hilbert space $\hilb$ with a \nsf\ trace $\tau.$
Let $\HH=\HH^*$ be affiliated with $\clN$ and let $\VV=\VV^* \in \LpN{1}.$ If $f \in \CPlus{2},$
then $f(\HH+\VV)-f(\HH) \in \LpN{1}$ and
$$
  \tau(f(\HH+\VV)-f(\HH)) = \int_\mbR f'(\lambda)\,d\Xi(\lambda),
$$
where the measure $\Xi$ is given by
$$
  \Xi(a,b) = \int_0^1 \tau(\VV E_{(a,b)}^{\HH_\rr})\,d\rr, \quad a,b \in \mbR.
$$
Here $\HH_\rr := \HH+\rr\VV, \ \rr \in [0,1]$ and $dE_\lambda^{\HH_\rr}$ is the spectral measure of $\HH_\rr.$
\end{thm}
Due to Lemma~\ref{L: trace of VE(l,s) is measurable} the measure $\Xi$ is well-defined.

\begin{proof}
  If $\phi(\lambda,\mu) = \alpha(\lambda)\beta(\mu),$ where $\alpha,\beta$ are continuous
bounded functions on $\mbR,$ then by the definition of the multiple operator integral
$$
  T^{\HH,\HH}_{\phi}(\VV) = \alpha(\HH)\VV\beta(\HH).
$$
Hence,
$$
  \tau\brackets{T^{\HH,\HH}_{\phi}(\VV)} = \tau(\alpha(\HH)\VV\beta(\HH))=\tau(\alpha(\HH)\beta(\HH)\VV).
$$
Since the function $\alpha(\cdot)\beta(\cdot)$ is bounded,
the simple spectral approximations to the bounded operator $\alpha(\HH)\beta(\HH)$
converge uniformly and so, after multiplying by $\VV$, converge in norm of $\LpN{1}.$
This implies that
$$
  \tau(\alpha(\HH)\beta(\HH)\VV) = \tau\brackets{\int_\mbR \alpha(\lambda)\beta(\lambda)\,dE_\lambda^\HH\, \VV}
  = \int_\mbR \alpha(\lambda)\beta(\lambda)\tau\brackets{dE_\lambda^\HH\, \VV}.
$$
Hence, for functions of the form $\phi(\lambda,\mu) = \alpha(\lambda)\beta(\mu),$
it follows that
\begin{gather}\label{F: (A)}
  \tau\brackets{T^{\HH,\HH}_{\phi}(\VV)} = \int_\mbR \phi(\lambda,\lambda)\tau\brackets{dE_\lambda^\HH\, \VV}.
\end{gather}

Let $(S,\Sigma,\nu)$ be a finite (complex) measure space, let $\alpha(\cdot,\cdot), \beta(\cdot,\cdot)$
be bounded continuous functions on $\mbR\times S$ and suppose that
$$
  \phi(\lambda,\mu) = \int_S \alpha(\lambda,\sigma)\beta(\lambda,\sigma)\,d\nu(\sigma) \quad \text{for all} \  (\lambda,\mu) \in \mbR^2
$$
is a representation of $\phi$ given by (\ref{F: BS representation}). Let $\phi_\sigma(\lambda,\mu) := \alpha(\lambda,\sigma)\beta(\mu,\sigma).$
It then follows from the definition of the multiple operator integral
that $T^{\HH,\HH}_\phi(\VV) = \int_S T^{\HH,\HH}_{\phi_\sigma}(\VV)\,d\nu(\sigma)$
and hence by Corollary~\ref{C: trace of MOI}
\begin{align*}
   \tau\brackets{T^{\HH,\HH}_{\phi}(\VV)}
   = \int_S \tau\brackets{T^{\HH,\HH}_{\phi_\sigma}(\VV)}\,d\nu(\sigma).
\end{align*}
It follows from (\ref{F: (A)}) that
\begin{align} \label{F: SSS}
   \tau\brackets{T^{\HH,\HH}_{\phi}(\VV)} & = \int_S \int_\mbR \phi_\sigma(\lambda,\lambda)\tau\brackets{dE_\lambda^\HH\, \VV}\,d\nu(\sigma) \notag\\
                 & = \int_\mbR \int_S \phi_\sigma(\lambda,\lambda)\,d\nu(\sigma)\tau\brackets{dE_\lambda^\HH\, \VV} \notag\\
                 & = \int_\mbR \phi(\lambda,\lambda)\tau\brackets{dE_\lambda^\HH\, \VV}.
\end{align}
The interchange of integrals in the second equality is justified by Lemma~\ref{L: the measure tau(AE(a,b)B) is count add}
and Fubini's theorem. 
Further, since $f \in \CPlus{2},$ it follows from Theorem~\ref{T: DalKreinIIb} applied to $\clE = \clL^1(\clN,\tau)$
that the Fr\'echet derivative $D_{\clL^1} f(\HH_\rr) = T^{\HH_\rr,\HH_\rr}_{\psif{f}}$ exists for all $r \in [0,1].$
By the continuity of the Fr\'echet derivative given by the estimate (\ref{F: norm estimate})
and the Newton-Leibnitz formula for the Fr\'echet derivative
(see e.g. \cite[Theorem 1.43]{Sch}) it follows that
\begin{gather*}
  \int_0^1 T^{\HH_\rr,\HH_\rr}_{\psif{f}}(\VV) \,d\rr = \int_0^1 D_{\clL^1} f(\HH_\rr)(\VV) \,d\rr = f(\HH+\VV) - f(\HH).
\end{gather*}
By Corollary~\ref{C: trace of MOI} and taking traces it follows that
\begin{gather}\label{F: NewtonLeibII}
  \int_0^1 \taubrs{T^{\HH_\rr,\HH_\rr}_{\psif{f}}(\VV)}\,d\rr = \tau(f(\HH+\VV) - f(\HH)).
\end{gather}

Since $\psif{f}$ is continuous, $\psif{f}(\lambda,\lambda)=f'(\lambda),$ so that (\ref{F: NewtonLeibII}) and (\ref{F: SSS}) imply
\begin{align*}
  \tau(f(\HH+\VV)-f(\HH)) & = \int_0^1 \int_\mbR \psif{f}(\lambda,\lambda)\tau\brackets{dE_\lambda^{\HH_\rr}\, \VV} \,d\rr \\
             & = \int_0^1 \int_\mbR f'(\lambda)\tau\brackets{dE_\lambda^{\HH_\rr}\, \VV} \,d\rr \\
             & = \int_\mbR f'(\lambda)\int_0^1 \tau\brackets{dE_\lambda^{\HH_\rr}\, \VV}\,d\rr,
\end{align*}
the interchange of the integrals in the last equality being justified
by Fubini's theorem \cite[VI.2]{Ja} via Lemma~\ref{L: the measure tau(AE(a,b)B) is count add}
and the fact that $f'$ is a bounded function.
\end{proof}
The next corollary in the case that $\clN = \clBH$ and $\tau=\Tr$ was established in \cite{BS72SM}.
\begin{cor}\label{C: The measure Xi is absolutely continuous} The measure $\Xi$ is absolutely continuous and the following equality holds
$$
  d\Xi(\lambda) = \xi(\lambda)\,d\lambda,
$$
where $\xi(t)$ is the spectral shift function for the pair $(\HH+\VV,\HH).$
\end{cor}
\begin{proof} From Theorem~\ref{T: Thm 3.1 ADS} and Theorem~\ref{T: BirmanSolII}, it follows that
$$
  \int_\mbR f'(\lambda)\,d\Xi(\lambda) = \int_\mbR f'(\lambda)\xi(\lambda)\,d\lambda
$$
for all $f \in C^\infty_\csupp(\mbR).$ Consequently, the measures $d\Xi(\lambda)$ and $\xi(\lambda)\,d\lambda$
have the same derivative in the sense of generalized functions. By \cite[Ch. I.2.6]{GSh}
there exists a constant $c$ such that
$$
  d\Xi(\lambda) - \xi(\lambda)\,d\lambda =c\cdot d\lambda.
$$
Since the measures $d\Xi(\lambda)$ and $\xi(\lambda)\,d\lambda$ are finite, it follows immediately, that $c = 0.$
\end{proof}

\section{Spectral shift and spectral flow}\label{S: SF and SSF}
\rndef{\HH}{D}
\rndef{\tHH}{D_1}

The second named author and \Phillips\ have established various
analytic formulae for spectral flow along a path of self-adjoint
unbounded Breuer-Fredholm operators affiliated with a semifinite von Neumann algebra.
For special choices of path suggested by the theory of the Krein
spectral shift function, one may study a spectral flow function on the real line:
$\mu\mapsto \sflow(\mu,\HH_0,\HH_1)$, $\mu\in\mbR$, where $\HH_1$ and $\HH_0$ differ
by a $\tau$-trace class operator and the function measures spectral flow from
$\HH_0-\mu$ to $\HH_1-\mu$. We now show that, under these
circumstances, the spectral flow function actually coincides with the Krein
spectral shift function.

Let us first recall preliminary material about spectral flow.
For more details see \cite{Ph96CMB,Ph97FIC,CP2}.
In these papers the notion of type II spectral flow is introduced
and an analytic approach is developed starting from ideas
of Getzler \cite{Ge93Top}. The new approach of these papers allows
the study of spectral flow between certain unbounded
Breuer-Fredholm operators affiliated with a general semifinite
von Neumann algebra \cite{CP2}.
We summarize the main features.

Let $\clN$ be a semifinite \vNa\ with a \nsf trace $\tau$ and $P,Q \in \clN$ be two infinite projections.
Let $\ker_Q T :=\ker T \cap Q(\clH).$
An operator $T \in P\clN Q$ is said to be $(P,Q)$ $\tau$-Fredholm if the subspaces
$\ker_Q T$ and $\ker_P T^*$ are $\tau$-finite and there exists a projection $P_1 \in \clN$ such that
$P_1 \le P,$ $\tau(P-P_1) < \infty$ and $P_1(\hilb) \subset T(\hilb).$ In this case $(P,Q)$-index of the operator $T$
is defined to be a number
$$
  \oind PQ(T) := \tau[\ker_Q T] - \tau[\ker_P T^*].
$$
Here $[\clK]$ denotes the projection onto the subspace $\clK \subseteq \hilb.$
If $P=Q=1$ we call $T$ just
$\tau$-Fredholm. For details see \cite{Br68MAnn,Br69MAnn,PR94JFA}.

Now, let $F \colon t\in [a,b] \mapsto F_t\in\clN$ be a norm continuous path of $\tau$-Fredholm operators
and $P_t = \frac 12(1+\sign(F_t)).$
If a partition $t_0=a < t_1 < \ldots < t_n = b$ of the segment $[a,b]$ is sufficiently small, then the operators
$P_{j-1}P_j \colon P_j\hilb \to P_{j-1}\hilb,$ $P_{j-1}P_j \in P_{j-1}\clN P_j,$ are $(P_{j-1},P_j)$ $\tau$-Fredholm
for $j = 1,\ldots,n$ ($P_j:=P_{t_j}$), so that the number
$$\sflow(\set{F_t}) := \sum\limits_{j=1}^n \oind {P_{j-1}}{P_j} (P_{j-1}P_j)
$$
is well-defined and does not depend on the partition $\set{t_j,j=1,\ldots,n}.$ Further,
if two paths $\set{F_t}$ and $\set{G_t}$ with the same ends points are norm homotopic,
then $\sflow(\set{F_t}) = \sflow(\set{G_t}),$ so that the
spectral flow $\sflow(F_0,F_1)$ depends only on the endpoints.

We recall the definition of a semifinite spectral triple (see e.g. \cite{CPRS1}).
\begin{defn} A semifinite spectral triple $\clAND$ is given by a Hilbert space $\hilb,$
a $*$-algebra $\clA \subset \clN$ where $\clN$ is a semifinite von Neumann algebra acting on $\hilb,$ and a densely
defined unbounded self-adjoint operator $\HH$ affiliated to $\clN$ such that \\
1) $[\HH,a]$ is densely defined and extends to a bounded operator for all $a\in \clA;$\\
2) $(\lambda -\HH)^{-1} \in \clKN$ for all $\lambda \notin \mbR,$ where $\clKN$ is the set of all $\tau$-compact operators.
\end{defn}
A spectral triple $\clAND$ is said to be $\theta$-summable if for all $t>0$ the operator $e^{-t\HH^2}$ is $\tau$-trace class.

Let $\clAND$ be a $\theta$-summable semifinite spectral triple, let $\VV \in \LpN{1}$ and let
$$
  \HH_r = \HH + r \VV, \quad r \in [0,1].
$$
The Carey-Phillips formula \cite{CP2} for spectral flow between $\HH_0=\HH$ and $\HH_1 = \HH+\VV$ is given by
\begin{multline*}
  \sflow(\HH_0,\tHH) = \epspi \int_0^1 \taubrs{\VV e^{-\eps{\HH_\rr}^2}}\,d\rr
     + \frac12 \brs{\eta_\eps(\HH_1) - \eta_\eps(\HH_0)}
   \\ + \frac12 \taubrs{[\ker(\HH_1)]-[\ker(\HH_0)]},
\end{multline*}
where the $\eta_\eps$-invariant of an unbounded self-adjoint operator $\HH,$ such
that $e^{-tD^2}$ is $\tau$-trace class for all $t>0,$ is defined
as \cite[Definition 8.1]{CP2}
$$
  \eta_\eps(\HH) := \frac 1{\sqrt \pi} \int_\eps^\infty \taubrs{\HH e^{-tD^2}} t^{-1/2}\,dt.
$$

Spectral flow may be interpreted as the `net amount' of spectrum crossing
zero while moving from $\HH_0$ to $\HH_1.$
So, it is natural to define the function $\sflow(\lambda; \HH_0, \HH_1):=\sflow(\HH_0-\lambda, \HH_1-\lambda)$
as spectral flow at a point $\lambda.$

It follows from the Carey-Phillips formula that
\begin{gather}\label{F: CP formula2}
  \sflow(\mu; \HH_0,\tHH) = \epspi \int_0^1 \taubrs{\VV e^{-\eps(\HH_\rr-\mu)^2}}\,d\rr
       + \frac12 \brs{\eta_\eps(\HH_1-\mu) - \eta_\eps(\HH_0-\mu)}
     \\ \notag +\frac12 \taubrs{[\ker(\HH_1-\mu)]-[\ker(\HH_0-\mu)]},   \quad \mu \in \mbR.
\end{gather}
The following theorem establishes a connection between the spectral shift function for the pair $(\HH_0+\VV,\HH_0)$
and the spectral flow function $\sflow(\,\cdot\,,\HH_0,\HH_0+\VV).$
\begin{thm}\label{T: SF thm} If $\VV \in \clN$ belongs to the $\tau$-trace class, then
$$
  \sflow(\mu; \HH_0, \HH_1) = \xi_{\HH_1, \HH_0}(\mu)+\frac12 \taubrs{[\ker(\HH_1-\mu)]-[\ker(\HH_0-\mu)]}
$$
for almost all $\mu \in \mbR.$
\end{thm}
\begin{proof} The spectral theorem and Lemma~\ref{L: net VA(alpha) converges to VA} implies
\begin{align*}
     \epspi \taubrs{\VV e^{-\eps(\HH_\rr-\mu)^2}}
      & = \epspi \taubrs{\VV\int_\mbR e^{-\eps(\lambda-\mu)^2}\,dE_\lambda^{\HH_\rr}} \\
      & = \int_\mbR j_\eps(\lambda-\mu)\,\taubrs{\VV dE_\lambda^{\HH_\rr}},
\end{align*}
where $j_\eps(x) = \epspi e^{-\eps x^2}, \ x \in \mbR.$
By Corollary~\ref{C: The measure Xi is absolutely continuous} and using the fact that the system $\set{j_\eps}$
is an approximate identity, we obtain that
\begin{multline*}
  \int_0^1 \int_\mbR j_\eps(\lambda-\mu)\,\taubrs{\VV dE_\lambda^{\HH_\rr}}\,d\rr
      = \int_\mbR j_\eps(\lambda-\mu)\frac d{d\lambda}\int_0^1 \taubrs{\VV E_\lambda^{\HH_\rr}}\,d\rr\,d\lambda
      \\ = j_\eps \ast \xi_{\HH_1,\HH_0}(\mu) \stackrel {L^1}\longrightarrow\xi_{\HH_1,\HH_0}(\mu) \quad \text{when} \quad \eps \to \infty.
\end{multline*}
The convergence in the last line can be justified by \cite[Chapter 3, section 5.6]{Rei}.

Since $\eta_\eps(\HH_j-\mu) \to 0, \ j = 0,1,$ when $\eps \to \infty,$
it follows from the Carey-Phillips formula (\ref{F: CP formula2}) that
$$
  \sflow(\mu; \HH_0, \HH_1) = \xi_{\HH_1,\HH_0}(\mu)+\frac12 \taubrs{[\ker(\HH_1-\mu)]-[\ker(\HH_0-\mu)]},   \quad \mu \in \mbR.
$$
\end{proof}

\mathsurround 0pt

\end{document}